\documentclass[3p]{elsarticle}

 \usepackage{graphics}
 \usepackage{graphicx}
 \usepackage{epsfig}

\usepackage{amssymb}
\usepackage{amsthm}
\usepackage{amsmath}

\usepackage{empheq}
\usepackage{epstopdf}
\usepackage{changepage}
\usepackage{rotating}
\usepackage{adjustbox}
\usepackage{color}
\usepackage{xcolor}
\usepackage{dsfont}
\usepackage{float}
\usepackage{lineno}
\usepackage{resizegather}
\usepackage{multicol}
\usepackage{booktabs} % mejora la calidad de las tablas y proporciona nuevas opciones
\usepackage{colortbl} % coloreado en tablas
\usepackage{multirow} % tablas con celdas de varias lineas
\usepackage{longtable} % tablas
\usepackage{hhline}   % mejora las lineas horizontales en las tablas
\usepackage{anyfontsize} % tamanos de fuente personalizados
\usepackage{everysel}
\usepackage{psfrag}
\usepackage{bbm}
\usepackage{bm}
\usepackage{pgf,tikz}
\usepackage[]{algorithm2e}
\usetikzlibrary{arrows}

\allowdisplaybreaks 

\newcommand{\R}{\mathbb{R}}

\newcommand{\x}{\chi}

\newcommand{\nstd}{\vec{n}_{j}}

\newcommand{\nv}{\vec{n}}

\newcommand{\TL}{\mathbf{L}}

\newcommand{\Gg}{ \bm{\mathcal{G}}}
\newcommand{\Mphi}{\bar{\bm{M}}}
\newcommand{\Mpsi}{\hat{\bm{M}}}
\newcommand{\Mphiref}{\bar{\mathcal{\bm{M}}}}
\newcommand{\Mpsiref}{\hat{\mathcal{\bm{M}}}}

\newcommand{\vv}{{\mathbf{v}}}
\newcommand{\xxi}{{\boldsymbol{\xi}}}

\newcommand{\xx}{ \mathbf{x}}
\newcommand{\TT}{\mathbf{T}}
\newcommand{\QQ}{\mathbf{R}}

\newcommand{\Ni}{ N_e}
\newcommand{\Ns}{ N_s}
\newcommand{\Nj}{ N_d}

\newcommand{\dx}{ d\xx}

\newcommand{\D}{\bm{\mathcal{D}}}
\newcommand{\Dref}{\mathcal{D}}
\newcommand{\U}{\bm{\mathcal{U}}}

\newcommand{\Q}{\bm{\mathcal{Q}}}

\newcommand{\diff}[2]{\frac{\partial {#1} }{\partial {#2} } }
\newcommand{\tphi}{\phi}
\newcommand{\tpsi}{\psi}
\newcommand{\tbphi}{\mathbf{\phi}}
\newcommand{\tbpsi}{\mathbf{\psi}}
\newcommand{\TTst}{\mathbf{T}}
\newcommand{\TTlst}{\TT_{\ell(j),j}}
\newcommand{\TTrst}{\TT_{r(j),j}}
\newcommand{\QQst}{\mathbf{R}}
\newcommand{\Nphist}{N_{\phi}}
\newcommand{\Npsist}{N_{\psi}}
\newcommand{\Nphi}{N_{\phi}}
\newcommand{\Npsi}{N_{\psi}}
\newcommand{\Gammast}{\Gamma}
\newcommand{\dxt}{d\xx}
\newcommand{\dSt}{ds}

\newcommand{\bbpi}{\bar{\mathbf{p}}}

\newcommand{\Qsm}{Q}
\newcommand{\Qsd}{Q}
\newcommand{\Qm}{\bar{Q}}
\newcommand{\Qd}{\hat{Q}}
\newcommand{\QQm}{\bar{\mathbf{Q}}}
\newcommand{\QQd}{\hat{\mathbf{Q}}}
\newcommand{\gradop}{\mathcal{G}}
\newcommand{\divop}{\mathcal{D}}

\newcommand{\Cm}{\bar{C}}

\newcommand{\gammad}{\hat{\gamma}}

\newcommand{\Lmd}{\mathcal{L}_{m \rightarrow d}}
\newcommand{\Ldm}{\mathcal{L}_{d \rightarrow m}}

\newcommand{\hbv}{\hat{\mathbf{v}}}

\definecolor{ttzzqq}{rgb}{0.2,0.6,0}
\definecolor{qqttcc}{rgb}{0,0.2,0.8}
\definecolor{qqttzz}{rgb}{0,0.2,0.6}
\definecolor{ffqqqq}{rgb}{1,0,0}
\definecolor{qqwuqq}{rgb}{0,0.39,0}
\definecolor{zzttqq}{rgb}{0.6,0.2,0}
\definecolor{qqqqff}{rgb}{0,0,1}
\definecolor{ttttqq}{rgb}{0.2,0.2,0}
\definecolor{qqwwtt}{rgb}{0,0.4,0.2}
\definecolor{ubqqys}{rgb}{0.29,0,0.51}
\definecolor{wwttqq}{rgb}{0.4,0.2,0}
\definecolor{uuuuuu}{rgb}{0.27,0.27,0.27}
\definecolor{qqzzff}{rgb}{0,0.6,1}
\definecolor{xdxdff}{rgb}{0.49,0.49,1}
\definecolor{ccwwqq}{rgb}{0.8,0.4,0}
\definecolor{ttqqqq}{rgb}{0.2,0,0}
\definecolor{qqzzcc}{rgb}{0,0.6,0.8}

\newcommand{\temp}{\theta}  				% temperature
\newcommand{\coeff}{\lambda} 				% diffusion coefficient
\newcommand{\dt}{\Delta t}					% time step
\newcommand{\Tm}{\bar{\temp}}               % temperature
\newcommand{\J}{\mathbf{J}}                 % Jacobian matrix

\newcommand{\UU}{\mathbf{q}}
\newcommand{\Um}{\bar{\UU}}
\newcommand{\Ud}{\hat{\UU}}
\newcommand{\MG}{\mathcal{M}_{\UU}}
\newcommand{\UUm}{\bar{\UU}}
\newcommand{\UUd}{\hat{\UU}}
\newcommand{\HH}{\bm{\mathcal{H}}}

\journal{Journal of Computational Physics}

\begin{document}

\begin{frontmatter}

\title{A new class of efficient high order semi-Lagrangian IMEX discontinuous Galerkin methods on staggered unstructured meshes}

\author[Unibz]{M. Tavelli\corref{cor1}}
\ead{mtavelli@unibz.it}
\author[usb,Fer]{W. Boscheri}
\ead{walter.boscheri@univ-smb.fr}
\cortext[cor1]{Corresponding author}
\address[Unibz]{Faculty of Engineering, Free University of Bolzano, Piazzetta Universit{\`a} 1, 39100 Bolzano, Italy}
\address[usb]{Laboratoire de Mathématiques UMR-5127 CNRS, Universit{\'e} Savoie Mont Blanc, 73376 Le Bourget du Lac, France}
\address[Fer]{Department of Mathematics and Computer Science, University of Ferrara, Italy}

% % % % % % % % % % % % % % % % % % % % % % % % % % % % % %
%                   Abstract                              %
% % % % % % % % % % % % % % % % % % % % % % % % % % % % % %
\begin{abstract}
In this paper we present a new high order semi-implicit discontinuous-Galerkin (DG) scheme on two-dimensional staggered triangular meshes applied to different nonlinear systems of hyperbolic conservation laws such as advection-diffusion models, incompressible Navier-Stokes equations and natural convection problems. While the temperature and pressure field are defined on a triangular main grid, the velocity field is defined on a quadrilateral edge-based staggered mesh. A semi-implicit time discretization is proposed, which separates slow and fast time scales by treating them explicitly and implicitly, respectively. The nonlinear convection terms are evolved explicitly using a semi-Lagrangian approach, whereas we consider an implicit discretization for the diffusion terms and the pressure contribution. High-order of accuracy in time is achieved using a new flexible and general framework of IMplicit-EXplicit (IMEX) Runge-Kutta schemes specifically designed to operate with semi-Lagrangian methods. To improve the efficiency in the computation of the DG divergence operator and the mass matrix, we propose to approximate the numerical solution with a less regular polynomial space on the edge-based mesh, which is defined on two sub-triangles that split the staggered quadrilateral elements. This allows for a fast computation of the DG operators and matrices without any need to store them for each mesh element. Due to the implicit treatment of the fast scale terms, the resulting numerical scheme is unconditionally stable for the considered governing equations because the semi-Lagrangian approach is the only explicit discretization for convection terms that does not need any stability restriction on the maximum admissible time step. Contrarily to a genuinely space-time discontinuous-Galerkin scheme, the IMEX discretization permits to preserve the symmetry and the positive semi-definiteness of the arising linear system for the pressure that can be solved at the aid of an efficient matrix-free implementation of the conjugate gradient method. We present several convergence results, including nonlinear transport and density currents, up to third order of accuracy in both space and time.
\end{abstract}

% % % % % % % % % % % % % % % % % % % % % % % % % % % % % %
%                   Keywords                              %
% % % % % % % % % % % % % % % % % % % % % % % % % % % % % %

\end{frontmatter}
 %\linenumbers

%                 Introduction                            %
% % % % % % % % % % % % % % % % % % % % % % % % % % % % % %
\section{Introduction}
Computational fluid mechanics plays an important role in environmental engineering applications that involve complex atmospheric flows, geophysical flows in oceans, rivers and lakes. Natural convection phenomena might also be relevant for instance in the simulation of incompressible flows in the atmosphere or even in the design of building isolation and solar collectors. These problems are mathematically modeled by the incompressible Navier-Stokes equations that take into account mass, energy and momentum conservation. Thermally driven phenomena can be described within this framework and, if slow flows and small density changes are assumed, the Boussinesq approximation can be applied, hence leading to the incompressible Navier-Stokes model coupled with an advection-diffusion equation for the temperature that determines the density fluctuations. 

A huge effort has been done to derive numerical methods able to solve small scale structures, see e.g. \cite{MNG15,MS18}. In the past decades a new class of high-order semi-implicit discontinuous Galerkin (DG) schemes was developed both on collocated \cite{GR10,Dol08,DF04,DFH07} and staggered meshes \cite{CL12,CCY13}, where the pressure contribution is taken implicitly while the convection terms are discretized explicitly, hence improving the overall efficiency of the scheme. In this approach, an implicit elliptic equation on the pressure is typically obtained, that results in a linear system which has to be solved for the unknown scalar pressure field. One advantage of the use of staggered (or dual) meshes in the DG framework is that, when applying integration by parts to the weak formulation of the governing equations, the staggered variables are continuous at the interface of the main grid and hence there is no need of a Riemann solver to compute the numerical fluxes. Consequently, the stencil size of the staggered scheme is smaller compared to classical collocated methods \cite{CS01,Cockburn98}. Edge-based staggering applied to fluid dynamics yields to the design of a compact discrete Laplace operator that involves only the direct neighbors of the main cell, see for example \cite{Casulli09,Casulli14,CG84,CS11}. Stability and efficiency of this family of semi-implicit methods on staggered meshes have been analyzed in \cite{CasulliCattani}, while an extension to mildly nonlinear pressure systems has been proposed in \cite{CZ09,CZ12}. High order DG schemes on staggered meshes were also employed for the solution of the shallow water and incompressible Navier-Stokes equations \cite{TD14,TD14sw}. There, the strategy adopted in finite-volume/finite-difference approximations \cite{Cas90,CC94} are modified and adapted to handle high order DG discretizations \cite{BR97,BO99,BO99b}, hence giving rise to a symmetric and semi-positive definite linear system for the pressure with arbitrary high order of accuracy in space. Applications of high order semi-implicit DG methods on Cartesian meshes can be found in \cite{DC13,FambriDumbser}. Those algorithms were also extended to arbitrary high order accuracy in space and time by using a genuinely space-time weak formulation \cite{TD16,TD17}. A drawback of the schemes presented in \cite{TD16,TD17} is that the introduction of time-dependent high order polynomials introduces an asymmetry in the linear system involving the pressure and the diffusion contributions. 

In most of the aforementioned cases, the nonlinear convective terms are discretized explicitly in order to avoid non-linearities in the implicit part. To treat those terms, a well-established explicit Godunov-type scheme is an attractive solution since it is simple to implement and stable under a Courant-Friedrichs-Lewy (CFL) condition. This restriction is based on the fluid velocity and thus it does not depend on the celerity, that is acceptable for several applications, see \cite{God59,Osher97,Munz1994}. An alternative way to solve the convective contribution is to use a so called semi-Lagrangian approach, where the convective terms are evaluated by following the Lagrangian trajectories of the flow over one time step and by subsequently interpolating the solution at the foot of the trajectory at the previous time level. It results in an explicit discretization that is unconditionally stable and therefore very appealing for practical and large scale applications \cite{Wel55,Wii59,ADERFSE,KW95,Yea09,VoronoiDivFree}. The semi-Lagrangian approach has been successfully employed for the convection terms in atmospheric flows \cite{Bonaventura2000} and more general incompressible fluids \cite{Bonaventura2018}, as well as for parabolic problems \cite{Bonaventura2014,BonaventuraFerretti2020}. Semi-Lagrangian approaches in the context of discontinuous Galerkin methods have been forwarded for incompressible flows as well as for weather prediction in \cite{RBS06,TBR13}. In a very recent work \cite{NatConv20}, a semi-implicit DG scheme was combined with a parallel semi-Lagrangian approach firstly designed in \cite{TB19}, and applied to natural convection problems. Second order of accuracy in time was obtained using a Crank-Nicholson time stepping technique while higher order time discretizations was achieved only by using a genuinely space-time DG formulation that results in a non symmetric linear system for the pressure and the temperature. Therefore, the efficient conjugate gradient method can not be used, thus resorting to a less efficient GMRES \cite{GMRES} solver.   

To treat problems with multiple time scales, as convection and pressure waves in fluid mechanics, an effective strategy is provided by the class of high order IMplicit-EXplicit (IMEX) methods, that have been investigated in the past decades \cite{AscRuuSpi,BP17,BP09,PR_IMEX,BPR2017}. 
%To treat problems with multiple time scales, as convection and pressure waves in fluid mechanics, an effective strategy is provided by the class of high order IMplicit-EXplicit (IMEX) methods, that have been investigated in the past decades \cite{AscRuuSpi,BP17,BP09,PR_IMEX,BPR2017}. 
In the IMEX approach, the slow and fast time scales are separated, hence permitting to devise efficient and asymptotic preserving integrators, meaning that the limit model of the governing equations is consistently reproduced at the discrete level \cite{JINAP1999, JP2001, KLARAP1999} when the stiffness parameter tends to zero, e.g. the Mach number in compressible flows or the Reynolds number. Very recently, IMEX schemes have been used with discontinuous Galerkin approximations for the solution of atmospheric flows and compressible fluids with non-conforming adaptive mesh refinement \cite{orlando2022_IMEXDG,orlando2023_IMEXDG}.

The IMEX schemes belong to the category of Method of Lines (MOL) time integrators, thus they naturally split the space and the time fluxes of the governing equations. As a consequence, once a first order in time semi-discrete scheme has been derived, the IMEX technique can be easily applied, independently from the specific space discretization that has been adopted. Indeed, IMEX time stepping schemes have been used with finite volumes \cite{BDLTV2020,BP21}, finite differences \cite{SIINS}, discontinuous Galerkin method \cite{BT22} or hybrid WENO-DG spatial operators \cite{BTC2023}. For the same reason, it is however not so easy and straightforward to couple IMEX time discretizations with a semi-Lagrangian approach due to the space-time nature of the Lagrangian trajectories that inevitably couple space and time. This problem was faced for the first time in a very recent paper \cite{TBP22} for both conservative and non-conservative semi-Lagrangian approaches. The applications are only limited to one-dimensional systems of PDE with finite volume and finite difference spatial discretizations.

In this paper we extend the non-conservative semi-Lagrangian IMEX (SL-IMEX) approach presented in \cite{TBP22} to two-dimensional unstructured staggered meshes and more challenging mathematical models. Our starting point is the semi-implicit DG method presented in \cite{NatConv20,TD14,TD14sw}, that is here modified as follows:
\begin{enumerate}
	\item we replace the genuinely space-time DG scheme with an IMEX time stepping technique combined with a DG discretization in space, that ultimately allows the symmetry of the pressure linear system to be retained, which is consequently solved relying on an efficient matrix-free conjugate gradient method;
	\item following \cite{BP21}, a semi-implicit IMEX scheme is devised, using the class of Runge-Kutta type integrators firstly introduced in \cite{BosFil2016}, hence improving the efficiency of the scheme by avoiding any need of fixed point Picard iterations as required in \cite{RDT20,TD14,TD16};
	\item the convective terms are discretized by means of the novel semi-Lagrangian methods forwarded in \cite{TBP22} for 1D governing equations, that are here extended to deal with DG spatial operators in multiple space dimensions on unstructured staggered meshes.
\end{enumerate}
Furthermore, one drawback of the methods presented in \cite{NatConv20,TD14,TD14sw} is that the definition of the basis functions on the edge-based staggered mesh does not allow a direct reconstruction of the differential operators with respect to a reference configuration. This makes the schemes very demanding in terms of memory requirements, since all the operators can be evaluated in the pre-processing stage but they are element dependent as underlined in \cite{TD16}. Along the lines of \cite{RDT20}, the regularity of the space spanned by the basis functions on the staggered control volumes can be reduced by splitting the dual cell into two sub-triangles. In this way the resulting functional space is continuous but with discontinuous derivative across the splitting edge. Thanks to this sub-grid strategy, all the main operators and matrices involved in the DG discretizations, namely divergence and mass matrix, will be computed on a limited set of  reference configurations to which the sub-triangles are mapped. We verify and validate the accuracy and the robustness of our approach with applications to a simple advection-diffusion equation, to the incompressible Navier-Stokes equations and to the natural convection model. 

The paper is organized as follows: the mathematical models are described in Section \ref{sec:governing_equations} while the staggered space discretization is detailed in Section \ref{sec:space_disc}. Section \ref{sec:SLIMEX} is devoted to introduce the semi-Lagrangian DG IMEX scheme, and we provide numerical results in Section \ref{sec:num_exp}. Finally, Section \ref{sec:conclusions} contains some conclusions and final remarks.

% % % % % % % % % % % % % % % % % % % % % % % % % % % % % %
%               Mathematical model                        %
% % % % % % % % % % % % % % % % % % % % % % % % % % % % % %
\section{Governing equations}\label{sec:governing_equations}
The new semi-Lagrangian IMEX discontinuous Galerkin schemes (SL-IMEX DG) are designed to operate on several governing equations, thus they exhibit great flexibility and a potentially large field of application. We consider a domain $\Omega \subset \mathds{R}^2$ with the generic position vector $\xx=(x,y)$, and $t \in \mathds{R}_0^+$ denotes the positive time coordinate.

Firstly, we present a simple advection-diffusion process of a scalar quantity $C=C(x,y,t)$ that moves according to a generic velocity field $\vv=(u,v)$:
\begin{equation}
	\diff{C}{t}+ \nabla \cdot \vv C=\nabla \cdot \left(\coeff \nabla C\right), \label{eqn:advdiff}
\end{equation}
where $\coeff$ represents a non-negative diffusion coefficient. 

Secondly, another well studied PDE system is provided by the incompressible Navier-Stokes equations that read
\begin{eqnarray}
	\nabla \cdot \mathbf{v}=0,\\
	\frac{\partial \mathbf{v}}{\partial t}+\nabla \cdot \mathbf{F_\vv} + \nabla p= \nabla \cdot \left( \nu \nabla \mathbf{v} \right),
\end{eqnarray}
where $p=P/\rho$ is the normalized pressure of the physical pressure $P$ with respect to the fluid density $\rho$, $\mathbf{F}_\vv=\vv \otimes \vv$ represents the tensor of nonlinear convective fluxes and $\nu=\mu/\rho$ is the kinematic viscosity coefficient. 

Lastly, we study an incompressible fluid where temperature/density gradients are taken into account. For small fluctuations of the temperature we can use the Boussinesq approximation \cite{NatConv20} instead of using a fully compressible model. In this case, the mathematical model is given by
\begin{eqnarray}
	\nabla \cdot \mathbf{v}=0, \label{eq:nc0.1}\\
	\frac{\partial \mathbf{v}}{\partial t}+\nabla \cdot \mathbf{F_\vv} + \nabla p= \nabla \cdot \left( \nu \nabla \mathbf{v} \right)+\left(  1- \beta \delta \temp\right) \mathbf{g}, \label{eq:nc0.2}  \\
	\diff{\temp}{t}+ \nabla \cdot \mathbf{F}_\temp =\nabla \cdot (\alpha \nabla \temp), \label{eq:nc0.3}
\end{eqnarray}
with $\alpha$ being the thermal diffusivity, that can be written in terms of the thermal conductivity $\kappa$, the heat capability $c_p$ and the density as $\alpha=\kappa/(\rho c_p)$. Moreover, $\beta$ is the thermal expansion coefficient of the fluid and $\mathbf{g}=(0,g)$ is the gravity acceleration vector, while $\temp$ represents the temperature and $\delta \temp$ is the temperature difference with respect to a reference value $\temp_0$, i.e. $\delta \temp=\temp-\temp_0$. The conservative fluxes in the temperature equation are given by the convective terms $\mathbf{F}_\temp = \vv \temp$.

% % % % % % % % % % % % % % % % % % % % % % % % % % % % % %
% % % % % % % % % % % % % % % % % % % % % % % % % % % % % %

\section{First order in time staggered semi-implicit DG scheme}\label{sec:space_disc}
The mathematical models presented in the previous section are numerically solved using a semi-implicit DG scheme on staggered unstructured meshes starting from the method forwarded in \cite{TD14,TD15,TD16}. More specifically, the unstructured staggered grid adopted here is the same such as in \cite{TD16} as well as the first order in time semi-implicit scheme. However, the nonlinear transport part is solved at the aid of the semi-Lagrangian schemes devised in \cite{TB19}, and a more efficient approach for the computation of the divergence operator and the mass matrix is proposed. In the following, we recall for completeness the notation and definitions adopted in this paper.

% % % % % % % % % % % % % % % % % % % % % % % % % % % % % %
% % % % % % % % % % % % % % % % % % % % % % % % % % % % % %
\subsection{Staggered unstructured mesh}\label{sec:mesh}
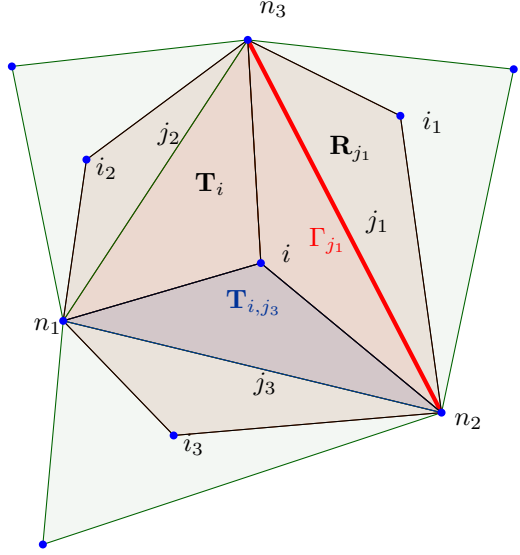
\begin{figure}[ht]
	\begin{center}
	    \begin{tikzpicture}[line cap=round,line join=round,>=triangle 45,x=0.6373937677053826cm,y=0.6177884615384613cm]
	    \clip(2.11,-8.53) rectangle (16.23,3.95);
	    \fill[color=zzttqq,fill=zzttqq,fill opacity=0.1] (5.19,-3.03) -- (9,3) -- (13,-5) -- cycle;
	    \fill[color=qqwuqq,fill=qqwuqq,fill opacity=0.05] (9,3) -- (14.49,2.37) -- (13,-5) -- cycle;
	    \fill[color=qqwuqq,fill=qqwuqq,fill opacity=0.05] (9,3) -- (4.13,2.43) -- (5.19,-3.03) -- cycle;
	    \fill[color=qqwuqq,fill=qqwuqq,fill opacity=0.05] (5.19,-3.03) -- (4.77,-7.83) -- (13,-5) -- cycle;
	    \fill[color=zzttqq,fill=zzttqq,fill opacity=0.1] (9.27,-1.79) -- (9,3) -- (5.67,0.43) -- (5.19,-3.03) -- cycle;
	    \fill[color=zzttqq,fill=zzttqq,fill opacity=0.1] (9.27,-1.79) -- (9,3) -- (12.15,1.37) -- (13,-5) -- cycle;
	    \fill[color=zzttqq,fill=zzttqq,fill opacity=0.1] (5.19,-3.03) -- (7.47,-5.49) -- (13,-5) -- (9.27,-1.79) -- cycle;
	    \fill[color=qqttzz,fill=qqttzz,fill opacity=0.1] (13,-5) -- (5.19,-3.03) -- (9.27,-1.79) -- cycle;
	    \draw [color=zzttqq] (5.19,-3.03)-- (9,3);
	    \draw [color=zzttqq] (9,3)-- (13,-5);
	    \draw [color=zzttqq] (13,-5)-- (5.19,-3.03);
	    \draw [color=qqwuqq] (9,3)-- (14.49,2.37);
	    \draw [color=qqwuqq] (14.49,2.37)-- (13,-5);
	    \draw [color=qqwuqq] (13,-5)-- (9,3);
	    \draw [color=qqwuqq] (9,3)-- (4.13,2.43);
	    \draw [color=qqwuqq] (4.13,2.43)-- (5.19,-3.03);
	    \draw [color=qqwuqq] (5.19,-3.03)-- (9,3);
	    \draw [color=qqwuqq] (5.19,-3.03)-- (4.77,-7.83);
	    \draw [color=qqwuqq] (4.77,-7.83)-- (13,-5);
	    \draw [color=qqwuqq] (13,-5)-- (5.19,-3.03);
	    \draw (9.51,-1.19) node[anchor=north west] {$i$};
	    \draw (12.41,1.67) node[anchor=north west] {$i_1$};
	    \draw (5.67,0.67) node[anchor=north west] {$i_2$};
	    \draw (7.47,-5.25) node[anchor=north west] {$i_3$};
	    \draw (11.25,-0.47) node[anchor=north west] {$j_1$};
	    \draw (6.95,1.43) node[anchor=north west] {$j_2$};
	    \draw (8.93,-3.93) node[anchor=north west] {$j_3$};
	    \draw (4.39,-2.79) node[anchor=north west] {$n_1$};
	    \draw (13.07,-4.83) node[anchor=north west] {$n_2$};
	    \draw (9.05,4.03) node[anchor=north west] {$n_3$};
	    \draw (7.71,0.35) node[anchor=north west] {$\TT_i$};
	    \draw [color=zzttqq] (9.27,-1.79)-- (9,3);
	    \draw [color=zzttqq] (9,3)-- (5.67,0.43);
	    \draw [color=zzttqq] (5.67,0.43)-- (5.19,-3.03);
	    \draw [color=zzttqq] (5.19,-3.03)-- (9.27,-1.79);
	    \draw [color=zzttqq] (9.27,-1.79)-- (9,3);
	    \draw [color=zzttqq] (9,3)-- (12.15,1.37);
	    \draw [color=zzttqq] (12.15,1.37)-- (13,-5);
	    \draw [color=zzttqq] (13,-5)-- (9.27,-1.79);
	    \draw [color=zzttqq] (5.19,-3.03)-- (7.47,-5.49);
	    \draw [color=zzttqq] (7.47,-5.49)-- (13,-5);
	    \draw [color=zzttqq] (13,-5)-- (9.27,-1.79);
	    \draw [color=zzttqq] (9.27,-1.79)-- (5.19,-3.03);
	    \draw (10.49,1.15) node[anchor=north west] {$\QQ_{j_1}$};
	    \draw [color=ffqqqq](10.07,-0.83) node[anchor=north west] {$\Gamma_{j_1}$};
	    \draw [line width=1.6pt,color=ffqqqq] (9,3)-- (13,-5);
	    \draw [color=qqttzz] (13,-5)-- (5.19,-3.03);
	    \draw [color=qqttzz] (5.19,-3.03)-- (9.27,-1.79);
	    \draw [color=qqttzz] (9.27,-1.79)-- (13,-5);
	    \draw [color=qqttzz](8.35,-2.17) node[anchor=north west] {$\TT_{i,j_3}$};
	    \draw (5.19,-3.03)-- (5.67,0.43);
	    \draw (5.67,0.43)-- (9,3);
	    \draw (9,3)-- (9.27,-1.79);
	    \draw (9.27,-1.79)-- (5.19,-3.03);
	    \draw (9,3)-- (12.15,1.37);
	    \draw (12.15,1.37)-- (13,-5);
	    \draw (13,-5)-- (9.27,-1.79);
	    \draw (13,-5)-- (7.47,-5.49);
	    \draw (7.47,-5.49)-- (5.19,-3.03);
	    \begin{scriptsize}
	    \fill [color=qqqqff] (5.19,-3.03) circle (1.5pt);
	    \fill [color=qqqqff] (9,3) circle (1.5pt);
	    \fill [color=qqqqff] (13,-5) circle (1.5pt);
	    \fill [color=qqqqff] (9.27,-1.79) circle (1.5pt);
	    \fill [color=qqqqff] (14.49,2.37) circle (1.5pt);
	    \fill [color=qqqqff] (4.13,2.43) circle (1.5pt);
	    \fill [color=qqqqff] (4.77,-7.83) circle (1.5pt);
	    \fill [color=qqqqff] (12.15,1.37) circle (1.5pt);
	    \fill [color=qqqqff] (5.67,0.43) circle (1.5pt);
	    \fill [color=qqqqff] (7.47,-5.49) circle (1.5pt);
	    \end{scriptsize}
	    \end{tikzpicture}
		\caption{Example of a triangular element $\TT_i$ composed by the edges $\{\Gamma_{j_1},\Gamma_{j_2},\Gamma_{j_3}\}$ with its neighbors. The quadrilateral dual element $\QQ_{j_1}$ is composed by the union of the left and right sub-triangles adjacent to the edge $\Gamma_{j_1}$, namely $\ell(j_1)=i$ and $r(j_1)=i_1$.}
		\label{fig.1}
	\end{center}
\end{figure}

The two-dimensional domain $\Omega \subset \mathbb{R}^2$ is paved using a set of non-overlapping triangles $\{\TT_i\}_{i=1}^{\Ni}$ where $\Ni$ is the total number of elements. Edges are indicated with $\Gamma_j$ for every $j\in [1,\Nj]$ with $\Nj$ being the total number of edges. For every edge which does not lie on the boundary, there exists a left and a right triangle that can be randomly assigned and they are denoted by $\ell(j)$ and $r(j)$, respectively. The standard orientation is intended from the left to the right element and a unit normal vector can be defined for every edge $j$ and is denoted with $\nv_j$. In a similar way, we define $S_i=\{j \in [1,\Nj] \,\, | \,\, \Gamma_j \mbox{ is an edge of } \TT_i\}$ as the set of edges that constitute the triangle $\TT_i$. Let $b_i\in \R^2$ be the center of mass of a triangle $\TT_i$, we define $\TT_{i,j}$ the triangle composed by the two nodes of $\Gamma_j$ and the point $b_i$. We then introduce an edge-based staggered mesh by defining $\QQ_j=\TT_{\ell(j),j}\cup \TT_{r(j),j}$ for every $j=[1,\Nj]$. To be noted that if $j$ refers to a boundary edge, then $\TT_{r(j),j}=\emptyset$ and then $\QQ_j=\TT_{\ell(j),j}$ is nothing but the left sub-triangle. However, $\QQ_j$ is a quadrilateral element if the edge $\Gamma_j$ is not a boundary edge. The notation is summarized in Figure $\ref{fig.1}$. From now on we refer to $\{\TT_i\}_{i=1}^{\Ni}$ as the \textit{main grid} and to $\{\QQ_j\}_{j=1}^{\Nj}$ as the \textit{dual grid}. The area of the main and dual control volume are labeled with $|\TT_i|$ and $|\QQ_j|$, respectively.

\subsection{Basis functions}\label{sec:basis_func}
\begin{figure}[!htbp]
	\begin{center}
		\includegraphics[width=0.60\textwidth]{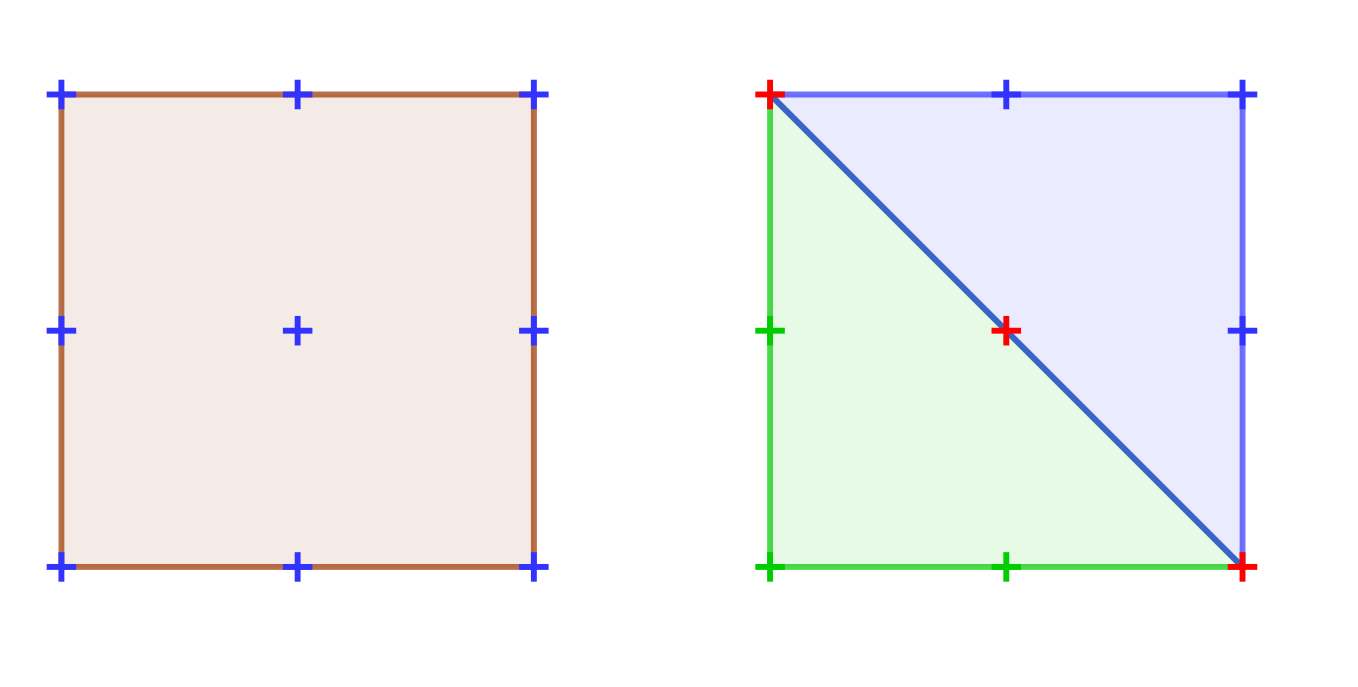}
		\caption{Example of functional space on the standard square $R_{ref}$ for $p=2$ with $9$ degrees of freedom (cross symbols). Tensor product of one dimensional basis functions (left) and reinterpretation in the context of agglomerated elements by split into two sub-triangles (right).}
		\label{fig.ADG}
	\end{center}
\end{figure}

On the main triangular mesh we can choose a polynomial basis up to a polynomial degree $p$ on a reference triangle $T_{ref}=\{(\xi,\eta) \in \R^{2} \,\, | \,\,  0 \leq \xi 
\leq 1, \,\,  0 \leq \eta \leq 1-\xi \}$ using a simple nodal approach, with $\xxi=(\xi,\eta)$ being the generic position vector in the reference space. Then, the connection from the reference to the physical element is performed using a linear mapping according to \cite{TD14,TD16}. We denote by $\{\phi_k \}_{k \in [1,N_\phi]}$ the set of basis functions on the reference space $T_{ref}$, where $N_\phi=\frac{(p+1)(p+2)}{2}$ is the size of the polynomial space. For the dual element, a reference basis on a unit square $R_{ref}$ can be built using tensor product of one-dimensional basis functions, see \cite{TD14}. However, this solution does not allow a direct connection between the physical and the reference discrete gradient/divergence operator because the mapping from the physical to the reference quadrilateral element is no longer linear and thus not analytically invertible, see \cite{TD14sw,TD14, TD16}. This means that all those operators become specific for each dual element with a significant demand in terms of memory, if they are computed and stored in a pre-processing stage. Therefore, we propose to relax the regularity of the solution representation in each dual element from being $C^\infty$ to $C^0$. This strategy has been already proposed in \cite{RDT20}, and it consists in separating the dual reference square into two sub-triangles that are simply the reference elements associated to $\TT_{\ell(j),j}$ and $\TT_{r(j),j}$ according to the orientation of the edge $\Gamma_j$. Figure \ref{fig.ADG} shows an example of the reference unit square $R_{ref}$ for $p=2$. As such, the dual basis functions can be expressed in terms of the main basis functions $\{\phi_k \}_{k \in [1,N_\phi]}$ defined on each reference sub-triangle. More in detail, the basis functions on the dual elements are labeled with $\{\psi_k \}_{k \in [1,N_\psi]}$ and defined for every $\xxi=(\xi, \eta) \in R_{ref}$ by
\begin{equation}
	\psi_{k}(\xxi)=\left\{
	\begin{array}{lc}
		\phi_{L(k,\xxi)}(\xxi) & \xi+\eta \leq 1 \\
		\phi_{R(k,\xxi)}(\xxi) & \xi+\eta > 1
	\end{array}
	\right. ,
\end{equation}
where $L(k,\xxi)$ and $R(k,\xxi)$ are two permutation functions for the indexes, see \cite{RDT20} for details. For example, let $p=1$, then $\psi(\xxi)=(\psi_1(\xxi),\psi_2(\xxi),\psi_3(\xxi),\psi_4(\xxi))$. The polynomial space becomes, in terms of $\phi$, $\psi(\xxi)=(\phi_1(\xxi),\phi_2(\xxi),\phi_3(\xxi),0)$ if $\xi+\eta\leq 1$ and $\psi(\xxi)=(0,\phi_2(1-\xxi),\phi_3(1-\xxi),\phi_1(1-\xxi))$ when $\xi+\eta> 1$.
To be noted that the nodal degrees of freedom that lie on the segment  $\xi+\eta=1$ are the same for both the left and the right triangle, thus ensuring $C^0$ continuity. The total number of degrees of freedom returns to be $N_\psi=(p+1)^2$, which is the same obtained in \cite{TD14}, but with a different interpretation that is based on a sub-triangulation of the reference dual control volume. The transformation from the reference space to the physical one and vice-versa is then driven by the mapping 
\begin{equation}
	\label{eqn.mapMain}
	T_i:\TT_i \rightarrow T_{ref} \qquad  \textnormal{and} \qquad T_i^{-1}:T_{ref} \rightarrow \TT_i
\end{equation}
on the main grid, which will be explicitly formulated in Section \ref{ssec.DGop} (see equation \eqref{eqn.linmap}). The map on the dual mesh $T_j$ is defined using two piecewise linear maps, hence obtaining the following mappings and their corresponding inverse: 
\begin{subequations}\label{eqn.mapDual}
	\begin{align}
	T_j(\TT_{\ell(j),j})\phantom{^{-1}}:\TT_{\ell(j),j}\rightarrow T_{ref} \qquad  \textnormal{and} \qquad  T_j(\TT_{r(j),j})\phantom{^{-1}}:\TT_{r(j),j}\rightarrow 1-T_{ref}, \label{eqn.MapDual1} \\
	T_j(\TT_{\ell(j),j})^{-1}:T_{ref} \rightarrow \TT_{\ell(j),j} \qquad  \textnormal{and} \qquad T_j(\TT_{r(j),j})^{-1}:1-T_{ref} \rightarrow \TT_{r(j),j}. \label{eqn.MapDual2}	
	\end{align}
\end{subequations}
We underline that, contrarily to what happen in \cite{TD14,TD16} where the Jacobian matrix of the transformation was space dependent, the Jacobian matrix is now piecewise constant in the two sub-triangles. In other words, the solution on the dual mesh is described using the agglomeration of two main elements, where we use the same basis functions adopted for the main triangular grid. A similar approach has been recently developed for handling DG schemes on general polygonal meshes \cite{ADERAFEDG}.

% % % % % % % % % % % % % % % % % % % % % % % % % % % % % %
% % % % % % % % % % % % % % % % % % % % % % % % % % % % % %
\subsection{First order in time semi-implicit DG discretization}\label{sec:semi_imp_dg} 
The time coordinate is indicated with $t$ and it is defined within the interval $[t_0,t_{end}]$, being $t_0$ and $t_{end}$ the initial and final time of the simulation, respectively. The time interval is discretized using $N+1$ points such that $t_0<t_1<\ldots < t_N=t_{end}$ with the generic time step given by $\Delta t_n=t^n-t^{n-1}$. Since the new SL IMEX-DG scheme is unconditionally stable, no CFL-type stability condition is required, and the time step is simply chosen according to the precision dictated by the physical phenomena and not by numerical restrictions.

Once the basis space for both the main and the dual grid is defined, we have to identify the variables that are assigned to the main and the dual grid as well. From now on, we adopt the convention that variables with a bar refer to the main grid while we use the hat symbol to indicate the dual variables. For a general quantity $\Qm$ defined on the main grid, i.e. $\Qm_i(\xx,t)=Q(\xx,t)|_{\TT_i}$, the numerical solution is written in terms of the polynomial basis $\{\phi_k\}_k$ as
\begin{eqnarray}
	\Qsm_i(\xx,t)      & = & \sum\limits_{l=1}^{\Nphi} \tphi_l^{(i)}(\xx) \Qm_{l,i}(t)=:\tbphi^{(i)}(\xx)\QQm_i(t), \qquad \forall i=1,\ldots \Ni, \label{eq:QmBFa}
\end{eqnarray}  
with the vector $\QQm_i=(\Qm_{1,i},\ldots , \Qm_{\Nphi,j})$ collecting all the degrees of freedom. The numerical solution of a variable defined on the dual grid $\Qd_j(\xx,t)=Q(\xx,t)|_{\QQ_j}$, is expressed accordingly by
\begin{eqnarray}
	\Qsd_j(\xx,t)      & = & \sum\limits_{l=1}^{\Npsi} \tpsi_l^{(j)}(\xx)
	\Qd_{l,j}(t)=:\tbpsi^{(j)}(\xx)\QQd_j(t), \qquad \forall j=1,\ldots \Nj, \label{eq:QdBFa}
\end{eqnarray} 
with the vector $\QQd_j=(\Qd_{1,j},\ldots , \Qd_{\Npsi,j})$ containing all the expansion coefficients. 
%Here the variables that identify the degrees of freedom on the main grid are represented with a bar while at the dual level they are identified with a hat.
Furthermore, the basis functions in the physical space $\phi_l^{(i)}$ and $\psi_l^{(j)}$ can be easily computed starting from the basis functions in the reference space and the considered physical element $\TT_i$ as $\phi_l^{(i)}(\xx)=T_i^{-1}(\phi_l(\xxi))$ and $\psi_l^{(j)}(\xx)=T_j^{-1}(\psi_l(\xxi))$, where the inverse mappings are provided in \eqref{eqn.mapMain} and \eqref{eqn.mapDual}, respectively.

Now, let us derive a first order in time semi-implicit scheme relying on the backward Euler time integrator, that will be generalized to the high order IMEX framework in Section \ref{sec:SLIMEX}. Moreover, a semi-Lagrangian discretization is employed for the convective terms along the lines of \cite{TB19}. 

We start by rewriting system \eqref{eqn:advdiff} as follows:
\begin{eqnarray}
	\frac{dC}{dt} &=&\nabla \cdot \gamma, \label{eqn:ad1.1} \\
	\gamma&=& \coeff \nabla C, \label{eqn:ad1.2}
\end{eqnarray}
where $dC/dt=\partial C /\partial t + \nabla \cdot \vv C$ is the total or Lagrangian derivative and $\gamma(\xx,t)$ is an auxiliary variable that represents the gradient of $C=C(\xx,t)$. We assign $C$ to the main grid while $\gamma$ is defined on the dual mesh. To obtain a weak formulation of the governing system, we multiply the first equation \eqref{eqn:ad1.1} by a test function $\phi_k^{(i)}$ and we integrate over the control volume $\TT_i$. Similarly, we multiply the second equation \eqref{eqn:ad1.2} by a dual test function $\psi_k^{(j)}$ and we integrate over the dual control volume $\QQ_j$. Thus, we obtain
\begin{gather}
	\int\limits_{\TT_i} \tphi_{k}^{(i)} \frac{dC}{dt} \, \dx=\int\limits_{\TT_i} \tphi_{k}^{(i)} \nabla \cdot \gamma \, \dx,\label{eqn:ad2.1} \\
	\int\limits_{\QQ_j} \tpsi_{k}^{(j)} \gamma \, \dx  = \int\limits_{\QQ_j} \tpsi_{k}^{(j)} \coeff \nabla C \, \dx. \label{eqn:ad2.2}
\end{gather}
Integrating by parts the term on the right hand side of equation \eqref{eqn:ad2.1} leads to
\begin{equation}
	\int\limits_{\TT_i} \tphi_{k}^{(i)} \frac{dC}{dt} \, \dx = \oint\limits_{\partial \TT_i}\tphi_{k}^{(i)} \gamma \cdot \nv_{i} \, \dSt -\int\limits_{\TT_i}\nabla \tphi_{k}^{(i)} \cdot \gamma \, \dxt ,
\end{equation}
where $\nv_{i}$ is the unit normal vector pointing out from the element boundary $\partial \TT_i$.
Note that since $\gamma$ is defined on the dual control volume, it is not continuous inside $\TT_i$. We then split the contribution of $\gamma$ within every sub-triangle $\TT_{i,j}$ with $j \in S_i$ and we use the definitions \eqref{eq:QmBFa} and \eqref{eq:QdBFa} to get the following first order in time semi-implicit scheme:
\begin{eqnarray}
	\int\limits_{\TT_i} \tphi_{k}^{(i)} \tphi_{l}^{(i)} \frac{\Cm_{l,i}^{n+1}-\Cm_{l,F}}{\Delta t}\dx&=&\sum\limits_{j \in S_i}\left(\int\limits_{\Gamma_j} \tphi_{k}^{(i)} \gamma \cdot \nv_{i,j} \, \dSt -\int\limits_{\TT_{i,j}}  \nabla \tphi_{k}^{(i)}  \cdot \gamma \, \dx\right)^{n+1} \nonumber \\ 
	&=& \sum\limits_{j \in S_i}\left(\int\limits_{\Gamma_j} \tphi_{k}^{(i)} \tpsi_{l}^{(j)} \nv_{i,j} \, \dSt -\int\limits_{\TT_{i,j}}  \nabla \tphi_{k}^{(i)}  \tpsi_{l}^{(j)} \dx\right)^{n+1} \cdot \gammad_{l,j}. \label{eqn:ad3.1}
\end{eqnarray}
Here, $\Cm_{l,F}=\Cm_{l,F}(\Cm_l^n,\vv)$ represents the value of $\Cm_l$ at the foot of the Lagrangian trajectory that has been obtained with the semi-Lagrangian scheme described later in Section \ref{sec:semilagr}, while $\nv_{i,j}=\nv_{i}|_{\Gamma_j}$ is nothing but the unit normal vector restricted to edge $\Gamma_j$, and it can be related to the standard normal vector as $\nv_{i,j}=\sigma_{i,j}\nv_j$ where $\sigma_{i,j}$ is a sign function that is $1$ if $i=\ell(j)$ and $-1$ otherwise:
\begin{equation} \label{eqn.sigma}
	\sigma_{i,j}=\frac{r(j)-2i + \ell(j)}{r(j)-\ell(j)}.
\end{equation}
The previous equation \eqref{eqn:ad3.1} is then written in a compact vector form as
\begin{equation}
	\Mphi_i \Cm_i^{n+1}=\Mphi_i \Cm^{n}_F+ \dt \sum\limits_{j \in S_i}\D_{i,j}\gammad_j^{n+1} \label{eq:ad4.1}
\end{equation}
where the mass matrix $\Mphi_i$ and the divergence operator $\D_{i,j}$ are defined component-wise as
\begin{eqnarray}
	\Mphi_i=\Mphi_{i;kl} &=& \int\limits_{\TT_i} \tphi_{k}^{(i)} \tphi_{l}^{(i)} \, \dx, \label{eq:Mi}\\
	\D_{i,j}=\D_{i,j;kl} &=& \int\limits_{\Gamma_j} \tphi_{k}^{(i)} \tpsi_{l}^{(j)} \cdot \nv_{i,j} \, \dSt -\int\limits_{\TT_{i,j}}  \nabla \tphi_{k}^{(i)}   \tpsi_{l}^{(j)} \, \dx. \label{eq:D}
\end{eqnarray}
Next, let us handle equation \eqref{eqn:ad2.2}. We do not perform integration by parts, but we point out that $C$ is discontinuous on $\Gamma_j$, therefore we propose to interpolate the gradient of $C$ in the sense of distributions across the edge $\Gamma_j$, thus obtaining
\begin{eqnarray}
	\int\limits_{\QQ_j} \tpsi_{k}^{(j)} \gamma_j^{n+1} \, \dx &=& \coeff \, \left( \int\limits_{\TTlst} \tpsi_k^{(j)} \nabla C_{\ell(j)} \, \dxt 
	+ \int\limits_{\TTrst}\tpsi_k^{(j)} \nabla C_{r(j)} \, \dxt 
	+\int\limits_{\Gammast_{j}}{\tphi_k^{(j)} \left(C_{r(j)}-C_{\ell(j)}\right) \nstd \, \dSt} \right)^{n+1} \nonumber \\
	\int\limits_{\QQ_j} \tpsi_{k}^{(j)} \tpsi_{l}^{(j)} \, \dx \, \gammad_{l,j}^{n+1} &=& \coeff \,  \left(\int\limits_{\TTlst} \tpsi_k^{(j)} \nabla \tphi_{l}^{(\ell(j))} \dxt  
	-\int\limits_{\Gammast_{j}} \tpsi_k^{(j)} \tphi_{l}^{(\ell(j))} \nstd \,\dSt\right) \, \, \Cm_{l,\ell(j)}^{n+1} \nonumber \\
	 &+&\coeff \, \left(\int\limits_{\TTrst}  \tpsi_k^{(j)} \nabla \tphi_{l}^{(r(j))} \, \dxt
	+\int\limits_{\Gammast_{j}}\tpsi_k^{(j)} \tphi_{l}^{(r(j))}    \nstd \, \dSt \right) \,  \, \Cm_{l,r(j)}^{n+1}.
\end{eqnarray}
The same formulation can also be obtained following the procedure originally introduced in \cite{BR97}. Even in this case, we rewrite the above equation in a compact form as
\begin{equation} 
	\Mpsi_j \gammad^{n+1}_j = \coeff \Q_{\ell(j),j}\Cm_{\ell(j)}^{n+1}+\coeff\Q_{r(j),j}\Cm_{r(j)}^{n+1}, \label{eq:ad4.2}
\end{equation}
where the mass matrix $\Mpsi_j$ and the gradient operator $\Q_{i,j}$ explicitly write
\begin{eqnarray}
	\Mpsi_j=\Mpsi_{j;kl}&=&\int\limits_{\QQ_j} \tpsi_{k}^{(j)} \tpsi_{l}^{(j)}\dx,  \label{eq:Mj} \\
	\Q_{i,j}=\Q_{i,j;kl} &=& \int\limits_{\TTst_{i,j}}\tpsi_k^{(j)} \nabla \tphi_{l}^{(i)}  \, \dxt-\int\limits_{\Gammast_j} \tpsi_k^{(j)} \tphi_{l}^{(i)}\sigma_{i,j} \nstd \, \dSt. \label{eq:ad5.2}
\end{eqnarray}
Let us underline that in the definition \eqref{eq:ad5.2} we have introduced the sign function \eqref{eqn.sigma} so that we obtain $\Q_{i,j}=\Q_{\ell(j),j}$ in \eqref{eq:ad4.2} if $\sigma_{i,j}=1$, i.e. if $i=\ell(j)$, and $\Q_{i,j}=\Q_{r(j),j}$ otherwise. In the semi-discrete scheme given by \eqref{eq:ad4.1} and \eqref{eq:ad4.2}, an implicit discretization is taken for $\hat{\gamma}^{n+1}$ and $\mathbf{\Cm}^{n+1}$, while the convective terms $\bar{\mathbf{C}}^{n}_F$ are explicitly computable relying on a semi-Lagrangian scheme (see Section \ref{sec:semilagr}). Substitution of the implicit $\gammad_j^{n+1}$ from \eqref{eq:ad4.2} into \eqref{eq:ad4.1} leads to
\begin{equation}
	\Mphi_i \Cm_i^{n+1}=\Mphi_i \Cm^{n}_F+\coeff \, \dt \, \sum\limits_{j \in S_i}\D_{i,j} \, \Mpsi_j^{-1} \, \left( \Q_{\ell(j),j}\Cm_{\ell(j)}^{n+1}+\Q_{r(j),j}\Cm_{r(j)}^{n+1}\right). \label{eq:ad5}
\end{equation}
This represents a linear system for the only unknown $\Cm_i^{n+1}$. Let us observe that the mass matrices $\Mphi$ and $\Mpsi$ are symmetric by construction, and that $\D=-\Q^\top$ (see equations \eqref{eq:D} and \eqref{eq:ad5.2}). Therefore, the linear system \eqref{eq:ad5} results to be symmetric and positive definite, see \cite{TD15,TD16} for further details.

For the Navier-Stokes case, we can define the pressure on the main grid and the velocity field on the dual one. Following the same reasoning we thus obtain 
\begin{gather}
	\sum_{j \in S_i}\D_{i,j}\hbv_j^{n+1}=0 \label{eq:ins.1},\\
	\Mpsi_{j} \TL_{h}^{\mathbf{v}}(\hbv_{j}) + \dt \, \left( \Q_{r(j),j} \bbpi_{r(j)}^{n+1}  + \Q_{\ell(j),j} \bbpi_{\ell(j)}^{n+1} \right) = 0, \label{eq:ins.2}
\end{gather}
where $\D$ and $\Q$ are the same operators defined in \eqref{eq:D} and \eqref{eq:ad5.2}, respectively, while  $\TL_{h}^{\mathbf{v}}(\hbv_{j})$ denotes the advection-diffusion operator which is computed here using a semi-Lagrangian approach for the transport part and an implicit solver for the diffusion fluxes, see \cite{NatConv20}. Following \cite{TD16}, a fractional step method is employed to solve the convective-diffusion subsystem on the main grid that accounts for three main steps: i) we first project the velocity field on the main grid using a high order $L_2$ projection; ii) then we solve the advection-diffusion equation \eqref{eq:ad5} for every component of the velocity vector; iii) the obtained solution is finally mapped back to the dual grid via a $L_2$ projection operator.

The model for the natural convection is also treated in a similar way, hence yielding the following scheme:
\begin{gather}
	\sum_{j \in S_i}\D_{i,j}\hbv_j^{n+1}=0 \label{eq:nc.1},\\
	\Mpsi_{j} \TL_{h}^{\mathbf{v}}(\hbv_{j}) + \dt \, \left( \Q_{r(j),j} \bbpi_{r(j)}^{n+1}  + \Q_{\ell(j),j} \bbpi_{\ell(j)}^{n+1} \right) = \dt \, \left( \Gg_j - \beta g \U_{r(j),j}\, \delta\Tm_{r(j)}^{n+1} - \beta g \U_{\ell(j),j}\, \delta\Tm_{\ell(j)}^{n+1} \right), \label{eq:nc.2} \\
	\Mphi_{i} \TL_{h}^{\temp}(\Tm_i, \bar{\vv})=0, \label{eq:nc.3} 
\end{gather}
with $\TL_{h}^{\temp}(\Tm_i, \bar{\vv})$ accounting for the advection-diffusion operator. In equation \eqref{eq:nc.2}, the integration of the source term leads to other operators that are given by
\begin{equation}
	\Gg_j=\Gg_{j;k}=\int_{\QQst_j}\tpsi_k^{(j)} \mathbf{g} \, \dxt, \label{eq:G}
\end{equation}
and
\begin{equation}
	\U_{j}=\U_{j;kl}=\int_{\TT_{i,j}}\tpsi_k^{(j)}  \tphi_{l}^{(i)}  \, \dxt. \label{eq:u}
\end{equation}
%that is defined when $j \in S_i$.
Let us notice that the same mathematical structure is exhibited for the advection-diffusion terms in the Navier-Stokes and the natural convection model, which is therefore numerically solved in the same manner by the operators $\TL_{h}^{\mathbf{v}}(\hbv_{j})$ and $\TL_{h}^{\temp}(\Tm_i, \bar{\vv})$, respectively. More in detail, the pure convection is computed relying on a semi-Lagrangian method, while the diffusion fluxes are implicitly solved. Indeed, this means that we are actually solving the same advection-diffusion system \eqref{eqn:advdiff} with the scheme given by \eqref{eq:ad4.1} and \eqref{eq:ad4.2}, applied to the velocity components for the Navier-Stokes system, and to the temperature for the natural convection model. The two operators $ \TL_{h}^{\mathbf{v}}(\hbv_{j})$ and $\TL_{h}^{\temp}(\Tm_i, \bar{\vv})$ can be written as $ \TL_{h}^{\mathbf{v}}(\hbv_{j})=\hbv_j^{n+1} - \mathbf{F}_{t}^{\mathbf{v}}\left( \bar{\vv}^{n}\right)$ and $\TL_{h}^{\temp}(\Tm_i, \bar{\vv})=\Tm_i^{n+1}-\mathbf{F}_{t}^{\temp}\left( \Tm^n,\bar{\vv}^{n}\right)$, where $\mathbf{F}_{t}^{\temp}$ and $\mathbf{F}_{t}^{\mathbf{v}}$ are the solution of the advection-diffusion equation in the form of \eqref{eq:ad5}, as fully detailed in \cite{NatConv20}.

Since we are solving the incompressible Navier-Stokes equations, in the discrete natural convection system \eqref{eq:nc.1}-\eqref{eq:nc.3} the pressure does not have a real evolutionary equation, but it is the one that guarantees the divergence-free condition \eqref{eq:nc.1}. In order to obtain such a pressure field at time $t^{n+1}$, we formally substitute the momentum equation \eqref{eq:nc.2} into the divergence-free condition \eqref{eq:nc.1}, hence retrieving the following linear system 
\begin{equation}
	\sum\limits_{j \in S_i}\D_{i,j}\Mpsi_{j}^{-1}\left(  \Q_{r(j),j} \bbpi_{r(j)}^{n+1} + \Q_{\ell(j),j} \bbpi_{\ell(j)}^{n+1}  \right) 
	= 
	\sum\limits_{j \in S_i}\D_{i,j}\Mpsi_{j}^{-1}\left(  \Gg_j 
- \beta g \U_{r(j),j}\, \delta\Tm_{r(j)}^{n+1} - \beta g \U_{\ell(j),j}\, \delta\Tm_{\ell(j)}^{n+1} \right)
	+\frac{1}{\Delta t} \sum\limits_{j \in S_i}\D_{i,j}\mathbf{F}^{\mathbf{v}}_t\left( \bar{\vv}^{n}\right) . \label{eq:pressure_final} 
\end{equation} 
This system is in general symmetric and semi-positive definite \cite{TD15,TD17} and can be solved using a matrix free implementation of the conjugate gradient method in order to obtain the new pressure field $\bbpi^{n+1}$. Once this is done, the new velocity can readily be updated from the discrete velocity equation \eqref{eq:nc.2}, thus yielding $\hbv^{n+1}$.

\subsection{Space DG mass matrices and differential operators} \label{ssec.DGop}
The scheme \eqref{eq:nc.2}-\eqref{eq:nc.3}-\eqref{eq:pressure_final} formally looks identical to the method presented in \cite{NatConv20}. However, the involved operators are quite different since the dual basis functions are not computed in the same way, as previously explained in Section \ref{sec:basis_func}. More precisely, the numerical solution obtained in \cite{NatConv20} is a smooth function inside each dual element, while in this work we require only that the solution is continuous and derivable within each dual cell, but with a discontinuity in the derivative, hence only retaining $C^0$ regularity. In the sequel we describe the construction of the main operators involved in the proposed scheme, namely the mass matrices and the divergence.

\subsubsection{Mass matrices}
Let us define the inverse linear mapping $T^{-1}(\xxi,\TT)$ with $\xxi \in T_{ref}$ which maps the reference triangle into a general triangle $\TT$, characterized by the three vertexes $(X_0,Y_0), (X_1,Y_1)$ and $(X_2,Y_2)$. Likewise, $T(\xx,\TT)$ maps the physical cell to the reference triangle, see the definition \eqref{eqn.mapMain}. We recall that the position vector in the reference and in the physical space are given by $\xxi=(\xi,\eta)$ and $\xx=(x,y)$, respectively. The map $T^{-1}(\xxi, \TT):(\xi,\eta) \mapsto (x,y)$ is defined by
\begin{equation}
	\begin{array}{ccc}
		x&=&X_0+\xi (X_1-X_0)+\eta (X_2-X_0) \\
		y&=&Y_0+\xi (Y_1-Y_0)+\eta (Y_2-Y_0) 
	\end{array},
	\label{eqn.linmap}
\end{equation}
where the Jacobian matrix $\J(T^{-1})$ is given by
\begin{equation}
	\J(T^{-1})=\frac{\partial \xx}{\partial \xxi} = \left(
	\begin{array}{cc}
		X_1-X_0 &  X_2-X_0 \\
		Y_1-Y_0 &  Y_2-Y_0 \\
	\end{array} \right).
\end{equation}
Note that $\J(T^{-1})$ does not depend by $\xxi$. Furthermore, we assume that the determinant of the Jacobian matrix is always positive, i.e. $|\J(T^{-1})|>0$, meaning that the cell area must obviously be positive in order to guarantee a valid computational grid. We introduce the notation $T_i^{-1}(\xxi):=T^{-1}(\xxi,T_i)$ and $T_i(\xx):=T(\xx,T_i)$, that are simply the maps applied to the triangle $T_i$. As already explained in Section \ref{sec:basis_func}, the basis functions on the physical space are defined as $\phi_k^{(i)}(\xx)=\phi_k(T_i(\xx))$ and hence
\begin{equation} \label{eqn.mapbasis}
	\phi_k^{(i)}(T_i^{-1}(\xxi))=\phi_k(T_i(T_i^{-1}(\xxi)))=\phi_k(\xxi).
\end{equation}
Therefore, using the above relation and a change of variable formula, the mass matrix on the main grid can be computed for the element $T_i$ as 
\begin{equation}
	\label{eqn.MassMatrixM}
	\Mphi_{i;kl}=\int\limits_{T_{ref}}\phi_k(\xxi) \, \phi_l(\xxi) \, |\J(T_i^{-1})| \, d\xxi=A_i \, \Mphiref_{kl},
\end{equation}
where
\begin{equation}
	\Mphiref_{kl}=\int\limits_{T_{ref}}\phi_k(\xxi) \, \phi_l(\xxi) \, d\xxi,
\end{equation}
and $A_i=|\J(T_i^{-1})|=2|\TT_i|$ is related to the area of $\TT_i$.
Similarly, we compute the dual mass matrix by splitting the quadrilateral cell into two sub-triangles, hence obtaining
\begin{eqnarray}
	\Mpsi_{j;kl}&=& \int\limits_{\QQ_j}\psi_k^{(j)}(\xx) \, \psi_l^{(j)}(\xx) \, d\xx, \nonumber \\
	&=& \int\limits_{\TT_{\ell(j),j}}\psi_k^{(j)}(\xx) \, \psi_l^{(j)}(\xx) \, d\xx
	+\int\limits_{\TT_{r(j),j}}\psi_k^{(j)}(\xx) \, \psi_l^{(j)}(\xx) \, d\xx.
\end{eqnarray}
The map $T_j$ restricted to a sub-triangle $\TT_{i,j}$ is nothing but the linear map between $\TT_{i,j}$ and $T_{ref}$ if $i=\ell(j)$ or between $\TT_{i,j}$ and $1-T_{ref}$ if $i=r(j)$, see the definition \eqref{eqn.MapDual1}. In this way, the above integrals can be evaluated on the reference triangle as
\begin{eqnarray}
	\label{eqn.MassMatrixD}
	\Mpsi_{j;kl}&=& \int\limits_{T_{ref}}\psi_k(\xxi) \, \psi_l(\xxi) \, |\J(T_{\ell(j)j}^{-1})| \, d\xxi
	+\int\limits_{1-T_{ref}}\psi_k(\xxi) \, \psi_l(\xxi) \, |J(T_{r(j)j}^{-1})| \, d\xxi \nonumber \\
	&=& A_{\ell(j)j} \, \Mphiref^L_{kl}+A_{r(j)j} \, \Mphiref^R_{kl},
\end{eqnarray}
with the Jacobian determinant $A_{ij}=|\J(T_{ij}^{-1})|=2|\TT_{i,j}|$ of the sub-triangle $\TT_{i,j}$ and the dual mass matrices
\begin{eqnarray}
	\Mpsiref^L_{kl} &=& \int\limits_{T_{ref}}\psi_k(\xxi) \, \psi_l(\xxi) \,d\xxi, \\
	\Mpsiref^R_{kl} &=& \int\limits_{T_{ref}}\psi_k(1-\xxi) \, \psi_l(1-\xxi) \, d\xxi.
\end{eqnarray}

As a consequence, the mass matrices \eqref{eqn.MassMatrixM} and \eqref{eqn.MassMatrixD} on the main and the dual cell, respectively, are written in terms of standard mass matrices $\Mphiref,\Mpsiref^L, \Mpsiref^R$ and geometric quantities $\{A_i\}_i$ and $\{A_{i,j}\}_{i,j}$. The standard mass matrices are not element-dependent, thus they are computed only once on the reference triangle, even for the dual cells thanks to the splitting into two sub-triangles. This makes the scheme remarkably faster compared to the numerical method forwarded in \cite{NatConv20}. Only the geometric quantities $\{A_i\}_i$ and $\{A_{i,j}\}_{i,j}$ have to be evaluated and stored for each element, which however are much less demanding in terms of memory requirements compared to the full mass matrices.

\subsubsection{Divergence and gradient operators}
The divergence and gradient operators are more challenging to be defined because they involve elements that belong to the main and to the dual grid simultaneously. Let us consider a main cell $i\in[1,\Ni]$ and a dual cell $j \in S_i$. 

If we use the map $T_i$ for the main cell, the sub-triangle $\TT_{i,j}$ can be mapped into one of the three sub-triangles $(T_{ref}^1,T_{ref}^2,T_{ref}^3)$ according to the corresponding local edge $\Gamma_{j_1},\Gamma_{j_2},\Gamma_{j_3}$. The reference sub-triangles are defined following their counterpart in the physical space, thus they are delimited by the two nodes of one edge of the reference triangle and the center of mass in the reference space that is $T_i(b_i)$. On the other hand, the sub-triangle $\TT_{i,j}$ can also be mapped into the left or the right sub-triangle of the dual reference element $R_{std}$ using the mapping \eqref{eqn.MapDual1}. 

\begin{figure}[!htbp]
	\begin{center}
		\begin{tabular}{c} \includegraphics[trim={0.5cm 0.5cm 0.5cm 1cm},clip,width=0.7\textwidth]{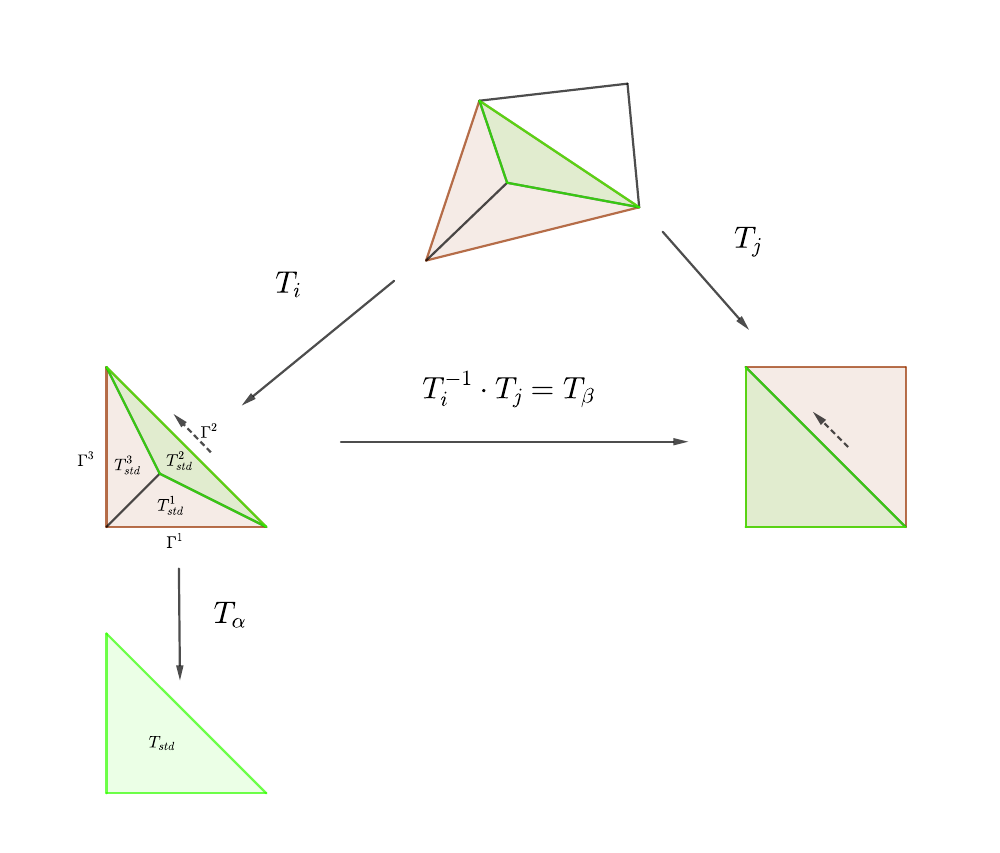} \\
		\end{tabular}
		\caption{Schematic of the reference maps. All the possible combinations are three for $T_\alpha$ and two for $T_\beta$}
		\label{fig:map2}
	\end{center}
\end{figure}

Let $T_\alpha:T_{ref}^\alpha\rightarrow T_{ref}$ be the linear transformation that maps the sub-triangle $T_{ref}^\alpha$ into $T_{ref}$ for $\alpha=1,2,3$ and $T_\beta=T_i^{-1}\circ T_j$. Furthermore, the three edges belonging to the reference triangle are naturally labeled with $\Gamma^{\alpha}$ for $\alpha=1,2,3$. A schematic view of the involved maps is shown in Figure \ref{fig:map2}.

Starting from the definition of the divergence operator $\D$ according to \eqref{eq:D}, we can identify a line and a volume contribution, hence $\D_{ij}=\D^{L}_{ij}-\D^{V}_{ij}$. By introducing the transformation $T_j^{E}:\Gamma_j \rightarrow [0,1]$, that maps the physical edge to the reference unit segment, we can explicitly write its inverse $T_j^{E,-1}(\tau,\Gamma):\tau \mapsto \xx(\tau)$ with $\tau \in [0,1]$ as
\begin{equation} \label{eqn.MapLine}
	 x(\tau)=X_1+\tau (X_2-X_1), \qquad y(\tau)=Y_1+\tau (Y_2-Y_1),
\end{equation}
where $(X_1,Y_1),(X_2,Y_2)$ are the coordinates of the two nodes defining the edge $\Gamma$. We also introduce the integration coordinates $s$ along the edge $\Gamma_j$ such that it runs from the first to the second node of the edge defined by the coordinates $(X_1,Y_1)$ and $(X_2,Y_2)$, respectively. The coordinate $s$ thus corresponds to the physical coordinates $\xx(s)$. Employing the edge map \eqref{eqn.MapLine}, we rewrite the line contribution of the divergence operator as 
\begin{eqnarray}
	\D^{L}_{ij}&=&\int\limits_{\Gamma_j}\phi_k^{(i)}(\xx(s)) \, \psi_l^{(j)}(\xx(s)) \, \vec{n}_j \sigma_{ij} \, dS \nonumber \\
	&=& |\Gamma_j| \, \vec{n}_j \int\limits_0^1 \phi_k^{(i)}\left( T_j^{E,-1}(\tau)\right)\psi_l^{(j)}\left( T_j^{E,-1}(\tau)\right)  \sigma_{ij} d\tau \nonumber \\
	&=& |\Gamma_j| \, \vec{n}_j \int\limits_0^1 \phi_k\left( T_i\left(T_j^{E,-1}(\tau)\right)\right)\psi_l\left( T_j\left(T_j^{E,-1}(\tau)\right)\right) \sigma_{ij} d\tau 
\end{eqnarray}
where we have used the definition of $\phi_k^{(i)}(\xx)=\phi(T_i(\xx))$ and $\psi_k^{(j)}(\xx)=\psi(T_j(\xx))$. Note that $ T_i\left(T_j^{E,-1}(\tau)\right)=\Gamma^\alpha$ for a proper $\alpha=1,2,3$. We can then define a new linear map $T_\alpha^{E,-1}(\tau):=T_i\left(T_j^{E,-1}(\tau)\right)$ so that $T_j\left(T_j^{E,-1}(\tau)\right)=T_j(T_i^{-1}(T_i(T_j^{E,-1}(\tau))))=T_\beta(T_\alpha^{E,-1}(\tau))$. This second transformation always maps the segment $[0,1]$ to the diagonal of $R_{ref}$ but with a different orientation, namely $T_j\left(T_j^{E,-1}(\tau)\right)=\Gamma_d$ if $i=\ell(j)$ or $T_j\left(T_j^{E,-1}(\tau)\right)=1-\Gamma_d$ otherwise. Here, $\Gamma_d$ is the diagonal of the reference dual element wich is always defined by the two nodes $(1,0)$ and $(0,1)$. The previous line integral then becomes
\begin{eqnarray}
	\D^{L}_{ij}&=& |\Gamma_j| \, \vec{n}_j \Dref^{L,\alpha,\beta}_{kl},
\end{eqnarray}
where
\begin{equation}
	\Dref^{L,\alpha,\beta}_{kl} = \int\limits_0^1 \phi_k\left(T_\alpha^{E,-1}(\tau)\right) \, \psi_l\left( T_\beta(T_\alpha^{E,-1}(\tau))\right) \, \sigma_{\beta} \, d\tau,
\end{equation}
that does not depend explicitly from $i$ and $j$ but has to be computed for $\alpha=1,2,3$ and $\beta=L,R$. Here, $\sigma_\beta=1$ if $\beta=L$ and $\sigma_\beta=-1$ otherwise.
%\bl{\begin{eqnarray}
%	\D^{L}_{ij}&=&\int\limits_{\Gamma_j}\phi_k^{(i)}(\xx(s)) \, \psi_l^{(j)}(\xx(s)) \, \vec{n}_j \sigma_{ij} \, dS \nonumber \\
%	&=&\int\limits_{\Gamma_j}\phi_k(T_i(\xx(s))) \, \psi_l(T_j(\xx(s))) \, \vec{n}_j \sigma_{ij} \, dS, \nonumber \\
%	&=&\int\limits_{\Gamma^\alpha} \phi_k(\xxi(s)) \, \psi_l(T_j(T_i^{-1}(\xxi(s)))) \, \vec{n}_j \sigma_{ij} \,|J(T_i^{-1})| ds, \nonumber \\
%	&=&\int\limits_0^1 \phi_k(T_\alpha^{-1}(s))\psi_l(T_\beta(T_\alpha^{-1}(s)))|J(T_i^{-1})||J(T_\alpha^{-1})|\vec{n}_j \sigma_{\beta} ds \nonumber \\
%	&=& |\Gamma_j| \vec{n}_j \Dref^{L,\alpha,\beta}_{kl},
%\end{eqnarray}
%where the reference line contribution is defined as
%\begin{equation}
%	\Dref^{L,\alpha,\beta}_{kl} = ...
%\end{equation}
%Note that the map to be adopted depends from the orientation of the element $\TT_{i,j}$ with respect to $\TT_i$ and $\QQ_j$. Hence, in order to obtain a global operator, all the combinations of $\alpha=1\ldots |S_i|$ and $\beta={L,R}$ has to be computed and stored in a pre-processing step. Accordingly, the value of $\sigma_\alpha$ depends from the orientation, namely  $\sigma_\alpha=1$ if $\beta=L$ and $\sigma_\alpha=1$ if $\beta=R$.} 
The volume contribution can be computed in a similar way as 
\begin{eqnarray}
	\D^{V}_{ij}&=&\int\limits_{\TT_{ij}}\nabla_\xx \phi_k^{(i)}(\xx) \, \psi_l^{(j)}(\xx) \, d\xx \nonumber \\
	&=& \int\limits_{T_{ref}^\alpha} \J(T_i)^\top \cdot\nabla_{\xxi} \,  \phi_k(T_i(T_i^{-1}(\xxi))) \, \psi_l(T_j(T_i^{-1}(\xxi))) \, |\J(T_{i}^{-1})| \,d\xxi \nonumber \\
	&=&A_{ij} \, \J_i \cdot \, \Dref^{V,\alpha,\beta}_{kl},
\end{eqnarray}
%\begin{eqnarray}
%	\D^{V}_{ij}&=&\int\limits_{\TT_{ij}}\nabla \phi_k^{(i)}(\xx) \, \psi_l^{(j)}(\xx) \, d\xx \nonumber \\
%	&=& \int\limits_{T_{ref}^\alpha} (\J(T_i))^{-1} \cdot\nabla_{\xxi} \,  \phi_k(T_i(T_i^{-1}(\xxi))) \, \psi_l(T_j(T_i^{-1}(\xxi))) \, |\J(T_{ij}^{-1})| \,d\xxi \nonumber \\
%	&=&A_{ij} \, (\J(T_i))^{-1} \, \Dref^{V,\alpha,\beta}_{kl},
%\end{eqnarray}
with the definition $\J_i=\J(T_i)^\top$ and the reference volume contribution 
\begin{eqnarray}
	\Dref^{V,\alpha,\beta}_{kl}=\int\limits_{T_{ref}}\nabla \phi_k(T_\alpha^{-1}(\xxi)) \, \psi_l(T_\beta(T_\alpha^{-1}(\xxi))) \, d\xxi	.
\end{eqnarray}
Combining back the volume and the line contribution we obtain the full divergence operator:
\begin{equation}
	\D_{ij}=|\Gamma_j|\, \vec{n}_j \Dref^{L,\alpha(i,j),\beta(i,j)}-A_{ij} \, \J_i \cdot \, \Dref^{V,\alpha(i,j),\beta(i,j)}.
\end{equation}
The following formulation of the gradient operator \eqref{eq:ad5.2} is derived with a similar reasoning, that is
\begin{equation}
	\Q_{ij}=|\Gamma_j| \vec{n}_j \Q^{L,\alpha(i,j),\beta(i,j)}-A_{ij} \,  \J_i \cdot \, \Q^{V,\alpha(i,j),\beta(i,j)},
\end{equation}
where we explicitly specify that $\alpha$ and $\beta$ are functions of $i,j$.
Note that also in this case we end up with a representation of the divergence and gradient operators in terms of $12$ standard operators defined in the reference system and some geometric quantities such as $ \J_i , \vec{n}_j$ and $|\Gamma_j|$ that are element-dependent.
\subsection{Semi-Lagrangian discretization of the convective terms} \label{sec:semilagr}
Let us consider the simple advection problem in non-conservative form:
\begin{equation} \label{eqn.adv}
	\diff{C}{t}+\vv \diff{C}{x}=0, \qquad C(\xx,0)= C_0(\xx). 
\end{equation}
By introducing the total or Lagrangian derivative $d/dt$, we have the identity 
\begin{equation}
	0=\frac{d}{dt}C(\xx,t)= \diff{C}{x}\vv+\diff{C}{t}, \label{eq:adv}
\end{equation}
meaning that the scalar field $C(\xx,t)$ is transported along the characteristic curves of the PDE \eqref{eqn.adv}. Therefore, 
the solution at time $t$ can be found following the characteristics backward in time up to $t=0$, that are governed by the trajectory equation
\begin{equation}
	\frac{d\xx}{dt}=\vv(\xx,t). \label{eq:lagrtrj}
\end{equation}
There are many ways to solve \eqref{eq:lagrtrj}, for example using high order trajectories \cite{TB19,VoronoiDivFree,ADERFSE} or simple piecewise linear paths.
%The previous equation \eqref{eq:lagrtrj} can be solved using high trajectories or piece-wise straight lines. 
Since the semi-implicit scheme is first order accurate in time so far, we use some simple straight lines to approximate the trajectories, which corresponds to use an Euler method for the solution of the ordinary differential equation (ODE) \eqref{eq:lagrtrj}. Once the backward time integration is carried out, we get the position $\xx_F$, that is typically referred to as the foot of the characteristic. There, the numerical solution has to be interpolated to finally obtain the contribution of the convective terms in the governing PDE, e.g. $\Cm^{n}_F$ in \eqref{eq:ad5}.

The semi-Lagrangian scheme is thus composed of two main steps: i) solution of the trajectory equation \eqref{eq:lagrtrj} and determination of the point $\xx_F$; ii) interpolation of the numerical solution at $\xx_F$, which requires the knowledge of the cell that contains the foot of the characteristic, hence implying a searching algorithm over the entire mesh. The second step can introduce a non negligible computational load, thus we pursue a different and very efficient strategy firstly developed in \cite{TB19}. Instead of computing the final point $\xx_F$ and performing a search within the elements of the whole mesh, we intersect the trajectory element by element, by solving the interpolation problem between the trajectory and the boundary of the cells. This permits to avoid the searching procedure because only the direct neighbors of the cell in which the trajectory has currently arrived are considered. % The time step for the backward time integration of the characteristics must therefore be chosen in such a way that no more than one cell per time step is crossed by the trajectory, thus we employ a sub-time stepping based on a classical CFL-type condition based on the transport velocity.
The time step for the backward time integration of the characteristics is driven by the intersection problem and guarantees that no more than one cell per time step is crossed by the trajectory, thus we employ a sub-time stepping based on the transport velocity.
  We underline that this is not a stability requirement of the scheme, but only a strategy to avoid a computationally very expensive global searching algorithm. Moreover, with this approach, the treatment of periodic boundaries are naturally included. 

Note also that, in the DG framework, equation \eqref{eq:adv} has to be evaluated for every degrees of freedom because of the weak formulation of the governing equations. In other words, one needs to compute the integral of the transported quantity at the foot of the characteristics and not just the centroid of the cell. This can be done relying on a Gauss quadrature for the evaluation of the integral. The Lagrange trajectory has then to be evaluated for every Gauss point, see \cite{TB19} for details.

\section{Semi-Lagrangian IMEX schemes}\label{sec:SLIMEX}
The introduction of a semi-Lagrangian discretization for the convective terms in combination with a staggered DG scheme was already proposed in \cite{NatConv20} using a Crank-Nicholson time discretization to achieve second order of accuracy in time. However, its combination with the family of Implicit-EXplicit (IMEX) Runge-Kutta schemes \cite{PR_IMEX,BosFil2016} is not trivial due to the space-time nature of the Lagrange trajectories, and, on the opposite side, the decoupling of space and time in Runge-Kutta integrators. This combination has been recently investigated in \cite{TBP22}, but only for one-dimensional problems. Here, we extend the algorithm presented in the aforementioned work to high order discontinuous Galerkin schemes on unstructured two-dimensional meshes.

In order to apply the IMEX framework to the natural convection model \eqref{eq:nc0.1}-\eqref{eq:nc0.3}, we first rewrite it in the form
\begin{equation}
	\diff{}{t}
	\left(
	\begin{array}{c}
		0 \\
		\vv \\
		\temp
	\end{array} 
	\right)
+\nabla \cdot
	\left(
\begin{array}{c}
	\vv \\
	\mathbf{F}_\vv -\nu \nabla \vv \\
	\mathbf{F}_\temp -\alpha \nabla \temp
\end{array} 
\right)
+
\nabla
\left(
\begin{array}{c}
	0 \\
	p \\
	0
\end{array} 
\right)
=
\left(
\begin{array}{c}
	0 \\
	(1-\beta \delta \temp)\mathbf{g} \\
	0
\end{array} 
\right),
\label{eq:imex_NC}
\end{equation}
which more compactly reduces to
\begin{equation}
	\partial_t (\boldsymbol{\tau} \, \UU) +\nabla \cdot \mathbf{f}(\UU, \nabla \UU) +\nabla \mathbf{h}(\UU)=\mathbf{S}(\UU), \label{eq:cform}
\end{equation}
where $\UU=(\vv, p,\temp)^\top$ is the state vector, $\mathbf{f}(\UU, \nabla \UU)$ contains the nonlinear convection and diffusion terms, $\nabla \mathbf{h}(\UU)$ accounts for the pressure gradients and $\mathbf{S}(\UU)$ represents the source terms in the momentum equation. In this framework, $\boldsymbol{\tau}=(0,1,1)^\top$ represents an indicator function that can assume the value $0$ or $1$, namely it is $\boldsymbol{\tau}^r=1$ if the component $\UU^r$ evolves in time, and zero otherwise. %The notation $\boldsymbol{\tau} \, \UU=\boldsymbol{\tau}^r \, \UU^r$ is then intended in a pointwise sense. 
Consequently, the first equation becomes the divergence-free constraint and it is not actually evolved but it is used to compute the fluid pressure that projects the numerical solution on the divergence-free manifold, according to the scheme \eqref{eq:pressure_final}.

Due to the adoption of a staggered mesh, we have two versions of $\UU$. Indeed, the state vector can be defined on the main or on the dual grid, hence yielding $\Um$ and $\Ud$, respectively. Let us introduce the subset $\MG=\{r \, | \, \UU^r \mbox{ is defined on the main grid}\}$, then the polynomial reconstructions given by \eqref{eq:QmBFa} and \eqref{eq:QdBFa} can be applied:
\begin{eqnarray}
	\Um_i^r(\xx,t)      & = & \sum\limits_{k=1}^{\Nphist} \tphi_k^{(i)}(\xx) \,
	\Um_{i;k}^r(t)=:\tbphi^{(i)}(\xx) \, \UUm_i^r(t) \label{eq:QmBF} \qquad \forall r\in \MG,  \\ 
	\Ud_j^r(\xx,t)      & = & \sum\limits_{l=1}^{\Npsist} \tpsi_k^{(i)}(\xx) \,
	\Ud_{i;k}^r(t)=:\tbpsi^{(i)}(\xx) \, \UUd_i^r(t) \label{eq:QdBF} \qquad \forall r\not\in \MG. \\
\end{eqnarray} 
We have also to provide two high order projection operators to transfer the solution from the dual to the main grid and vice-versa. They write
\begin{eqnarray}
	\Ldm(\UUd^r)	\qquad &:& \qquad \UUm_i^r = \Mphi^{-1}_i \sum_{j \in S_i} \U^\top_{i,j} \, \UUd_j^r \qquad \forall i=1 \ldots N_i, \label{eq:d2m}\\
	\Lmd(\UUm^r)	\qquad &:& \qquad \UUd_j^r = {\Mpsi}_j^{\, -1} \left( \U_{\ell(j),j} \, \UUm^r_{\ell(j)} + \U_{r(j),j} \, \UUm^r_{r(j)}\right) \qquad \forall j=1 \ldots N_j, \label{eq:m2d}
\end{eqnarray}
where $\U$ is the operator defined in \eqref{eq:u}, $\Mphi$ and $\Mpsi$ are the mass matrices given by \eqref{eq:Mi} and \eqref{eq:Mj}, respectively.
The definitions \eqref{eq:d2m} and  \eqref{eq:m2d} are valid for any scalar quantity formally defined on the dual and on the main grid, respectively. In order to apply the IMEX approach, we further cast the governing equations in the form of an autonomous system
\begin{eqnarray}
	\frac{d\UU}{dt}&=&\HH^P(\UU^E, \UU^I), \label{eq:IMEX1} \\
	\frac{d\xx}{dt} &=& \HH^L(\UU^E, \UU^I),  \label{eq:IMEX2}
\end{eqnarray}
with the introduction of the total derivative of $\UU$, as typically done in the Method Of Lines (MOL) framework. Moreover, the trajectory equation \eqref{eq:IMEX2} is also embedded in the IMEX procedure because it will provide the fluxes, i.e. the velocity, for the semi-Lagrangian scheme. Here, the index $I$ refers to implicit variables and $E$ is there to indicate explicit variables. The quantities $\HH^P,\HH^L$ represent the spatial approximation of the differential operators as well as the sources in the governing equations. Following the first order semi-implicit discretization \eqref{eq:nc.1}-\eqref{eq:nc.3}, we have
\begin{eqnarray}
	\HH^P=\left(
	\begin{array}{c} \HH^{P,\vv} \\ \HH^{P,\temp} \end{array}
	\right) = \left(
	\begin{array}{c}
		(\nabla \cdot \nu \nabla \theta)^I+\nabla p^I -((1-\beta \delta C)\mathbf{g})^E \\
		(\nabla \cdot \alpha \nabla \temp)^I
	\end{array}
	\right) \qquad \mbox{and} \qquad \HH^L=-\vv^E. \label{eqn.HsHL}
\end{eqnarray}
The pressure enters the fluxes as a constraint through the divergence-free condition that can be directly excluded in the system \eqref{eq:IMEX1}. Nevertheless, its contribution appears in $\HH^P$ and hence the pressure field needs to be computed at any stage.
%Note that the system does not include the evolution of the pressure  in \eqref{eq:IMEX1} since for incompressible Navier-Stokes $\gamma_r=0$.
This form follows directly from the one presented in \cite{BDLTV2020}, but here we are considering two staggered grids, implying that equation \eqref{eq:IMEX1} can be solved on the dual or on the main grid. Since the dual and main variables coexist on the same grids, either the dual or the main grid, the system \eqref{eq:IMEX1} is formulated for both meshes, so that the computation of the fluxes and the variables is performed only on the correct mesh where the evolution equation is formally defined (e.g. according to the location of $\UU$ on the grid), while the remaining parts are evaluated as high order projection with \eqref{eq:d2m} and \eqref{eq:m2d}:
\begin{eqnarray}
	\frac{d\Um^r}{dt}&=&\bar{\HH}^P(\UU^E, \UU^I) \qquad r\in \MG, \\
	\frac{d\Ud^r}{dt}&=&\hat{\HH}^P(\UU^E, \UU^I) \qquad r\notin \MG.
\end{eqnarray}
In the specific case of the natural convection model, the main and the dual contributions become
\begin{equation}
	\bar{\HH}^P=\left(
	\begin{array}{c}
		\Ldm(\hat{\HH}^{P,\vv}) \\
		(\nabla \cdot \alpha \nabla \temp)^I 
	\end{array}
	\right), \qquad 
		\hat{\HH}^P=\left(
	\begin{array}{c}
		(\nabla \cdot \nu \nabla \temp)^I-\nabla p^I +((1-\beta \delta \temp)\mathbf{g})^E \\
		\Lmd(\bar{\HH}^{P,\temp})
	\end{array}
	\right),
	\label{eqn.HPmd}
\end{equation}
and the two state vectors are given by
\begin{eqnarray}
	\Um=\left(
	\begin{array}{c}
		\Ldm(\hat{\vv}) \\
		\Tm
	\end{array}
	\right), \qquad 
	\Ud=\left(
	\begin{array}{c}
		\hat{\vv} \\
		\Lmd(\Tm) \label{Hsmd}
	\end{array}
	\right).
	\label{eq:statevecprj}
\end{eqnarray}
Note that since $\temp$ is formally defined on the main grid and $\vv$ is defined on the dual grid, in principle one would have to compute both terms $\bar{\HH}^P$ and $\hat{\HH}^P$. In practice, this is not the case because we only evaluate $\hat{\HH}^{P,\vv}$ and $\bar{\HH}^{P,\temp}$, then the remaining two contributions are simply obtained via projection, see equation \eqref{eqn.HPmd}.
% then $\Qm^C=C$ and $\Qd^\vv=\vv$, but we still compute $\Qm^\vv=\Ldm(\Qd^\vv)$ and $\Qd^C=\Lmd(\Qm^C)$.
%Let $\gradop$ and $\divop$ be the grad operators defined in \eqref{eq:ins.1} and \eqref{eq:ins.2}. 
%We remark that the gradient operator $\gradop$ applied to a scalar field on the main mesh returns the gradient at the dual level, see equation \eqref{eq:ad4.2}. On the contrary, the operator $\divop$ returns the divergence field on the main mesh, starting from a two-dimensional vector field defined at the dual level, see equation \eqref{eq:ins.1}.  
%%\begin{eqnarray}
%%	\gradop(\Qd^c) \qquad &:& \qquad \Mpsi_{j}^{-1}\left(  \Q_{r(j),j} \Qm^c_{r(j)} + \Q_{\ell(j),j} \Qm^c_{\ell(j)}  \right) \qquad \forall j=1 \ldots N_j , \\
%%	\divop(\Qm^c) \qquad &:& \qquad \sum\limits_{j \in S_i}\D_{i,j}\Qd^c_j \qquad \forall i=1 \ldots N_i , \\
%%\end{eqnarray}
%The Laplacian operator is then evaluated as $\mathcal{L}(\Um)=\divop(\gradop(\Um))$, yielding a symmetric operator from the main to the main grid, see e.g. \cite{TD15,TD17}. \bl{***** UP to here *****}
	
The scheme \eqref{eq:nc.1}-\eqref{eq:nc.3} written in the form of an autonomous system as \eqref{eq:IMEX1}-\eqref{eq:IMEX2} leads to a first order time marching algorithm. To achieve higher order of accuracy, we rely on IMEX Runge-Kutta schemes. These ODE integrators are typically described in terms of a double Butcher tableau of the form
	\begin{equation}
		\begin{array}{c|c}
			\tilde{c} & \tilde{A} \\ \hline & \tilde{b}^\top
		\end{array} \qquad
		\begin{array}{c|c}
			c & A \\ \hline & b^\top
		\end{array},
		\label{eqn.butcher}
	\end{equation}
where  $(\tilde{A},A) \in \R^{\Ns \times \Ns}$ and $(\tilde{c},c,\tilde{b},b) \in \R^{\Ns}$ depend on the number of stages $\Ns$. The left tableau represents the explicit scheme, thus $\tilde{A}=\tilde{a}_{ij}$ is a strictly lower triangular matrix, while the right tableau provides the coefficients for the implicit scheme with $A={a}_{ij}$, hence having non-zero entries on the diagonal. If ${a}_{N_s \, j}=b_j$ the IMEX scheme is Stiffly Accurate (SA), which is a remarkable property for the development of asymptotic preserving (AP) schemes \cite{BPR2017,BP21}, and the final numerical solution coincides with the last Runge-Kutta stage. In this work, the family of semi-implicit IMEX Runge-Kutta schemes is adopted, firstly introduced in \cite{BosFil2016}, and subsequently tested for multidimensional compressible and incompressible flows, see e.g. \cite{BFR16,BP21,TBP22}. We indicate with $R=0,1,2$ the degree of the time integrator, according to the specific schemes reported in \ref{app.imex}, implying that the DG discretizations can achieve up to third order of accuracy in time (for $R=2$).
	
Next, we apply semi-Lagrangian IMEX approach proposed in \cite{TBP22} for the discretization of the convective-viscous sub-system, that is solved employing a fractional step method. The implicit viscous sub-system has already been inserted in the contribution $\HH^P$ of \eqref{eqn.HsHL}, while for the explicit convective sub-system we have now to take into account the ODE fluxes $\HH^L$ introduced in \eqref{eqn.HsHL}, which will ultimately provide the foot of the characteristics where we need to interpolate the numerical solution. We remark that this scheme makes use of the velocity field projected on the main mesh, i.e. $\bar{\vv}(\xx)$, as suggested in \cite{TD16} and adopted in \cite{NatConv20}. For every IMEX stage $s=1 \ldots \Ns$, after initializing $\xx_1^E=\xx_1^I=\xx$, we compute
	\begin{eqnarray}
		\xx^E_{j+1}&=&\xx^E_{j}- \dt \, \tilde{a}_{s,j} \bar{\HH}^L_j(\xx) \qquad \forall j=1\ldots s-1, \\
		\xx^I_{j+1}&=&\xx^I_{j}- \dt \, a_{s,j} \bar{\HH}^L_j(\xx) \qquad \forall j=1\ldots s-1.
	\end{eqnarray}
Next, the explicit fluxes for the trajectory equation \eqref{eq:IMEX2} are evaluated at $\xx^E_s$, that is
	\begin{equation}
			\bar{H}^L_s=-\bar{\vv}(\xx^E_s).
	\end{equation}
	The last step of the implicit trajectory can then be readily updated:
	\begin{equation}
		\xx^I_{s}=\xx^I_{s} + \dt \, a_{s,s} \bar{H}^L_s(\xx).
	\end{equation}
The convective contribution for the stage $s$ is therefore provided by the interpolation of the numerical solution at $\xx^I_{s}$. This makes a strong coupling between space and time, which does not allow anymore to separate the purely spatial fluxes $\bar{H}^L_j(\xx)$ from the convective contribution, as pointed out in \cite{TBP22}. Therefore, the idea is to move also the spatial flux $\bar{H}^s_j(\xx)$ using the Lagrangian trajectories, so that it can be shifted with a time fractional step $w_s$ at the same time level of the advection contribution in the Runge-Kutta stages. To this purpose, we introduce the new shifted fluxes $\{\bar{\HH}^*\}_{j=1}^{s-1}$ defined as
	\begin{equation}
		\bar{\HH}^*_j:=\bar{\HH}^P_j(\xx+w_j\Delta t \bar{\HH}^L_j), \qquad \forall j=1\ldots s-1,
	\end{equation}
	where
	\begin{equation}
		w_j   = \sum\limits_{r=1}^j a_{j,r}-\sum\limits_{r=1}^{j} a_{s,r}.
	\end{equation}
Here, the notation $\bar{\HH}^P(\xx_F)$ is a compact way to describe the evaluation of $\bar{\HH}^P$ at the foot of the Lagrangian trajectory, i.e. at $\xx_F=\xx_F(\xx)$. Since all the involved quantities are approximated as polynomials in the DG framework, the shifted fluxes have to be computed in an integral sense and evaluated on every Gauss quadrature point $g=1\ldots N_g$ of coordinate $\xxi_g$ in the reference element $T_{ref}$ and weight $\omega_g$:
	\begin{eqnarray}
		\bar{\HH}^P_{i;k}(\xx_F)&=&\Mphi_{kl}^{-1}\int\limits_{\TT_i}\phi_l^{(i)}(\xx) \, \bar{\HH}^P(\xx_F(\xx)) \, d \xx \nonumber \\
		&=& \Mphi_{kl}^{-1}\sum\limits_{g=1}^{N_g}\phi_l^{(i)}(\xx(\xxi_g))\bar{\HH}^P(\xx^I(\xx(\xxi_g))) \, \omega_g  \, |\J(T^{-1})| \qquad \forall i=1\ldots \Ni, \quad k=1\ldots \Nphi. \label{eq:intstage}
	\end{eqnarray}
	Note that this is slightly different with respect to what presented in \cite{TBP22} since the fluxes are directly transported up to $\sum\limits_{r=1}^{j} a_{s,r}\Delta t$, whereas in \cite{TBP22} two projections were considered, namely up to $\sum\limits_{r=1}^{j-1} a_{s,r}\Delta t$ and then $a_{s,j}\Delta t$. This allows us to avoid one projection step and thus one re-integration stage. Once all the fluxes at the previous Runge-Kutta stage are computed and properly shifted, we can compute the intermediate solution accordingly for every stage $s=2 \ldots N_s$:
	\begin{eqnarray}
		\UU^{I}_s(\xx)&=&\UU^n(\xx)+\sum\limits_{j=1}^{s-1}\Delta t \, a_{s,j} \, \HH^*_j, \\
		\UU^{E}_s(\xx)&=&\UU^n(\xx)+\sum\limits_{j=1}^{s-1}\Delta t \, \tilde{a}_{s,j} \, \HH^*_j,
	\end{eqnarray}
    which hold true for both $(\bar{\UU}^{I}_s(\xx),\hat{\UU}^{I}_s(\xx))$ and $(\bar{\UU}^{E}_s(\xx),\hat{\UU}^{E}_s(\xx))$.
	%Recall that the computation of $\UU^n(\xx)$ requires an integration stage as described by the quadrature formula \eqref{eq:intstage}. 
	At this point, we need to evaluate the semi-implicit fluxes of the current Runge-Kutta stage. Therefore, we solve the linear system associated to the diffusion sub-system:
	\begin{eqnarray}
		\Tm^*-a_{s,s} \, \Delta t \, \nabla \cdot \alpha\nabla \Tm^* &=&\Tm^I_s,\\
		\bar{\vv}^*-a_{s,s} \, \Delta t \, \nabla \cdot \nu \nabla \bar{\vv}_s^* &=&\Ldm(\hat{\vv}^I_s).	
	\end{eqnarray}
	After computing $\vv^*=\Lmd(\bar{\vv}^*)$, we solve the linear system \eqref{eq:pressure_final} associated to the pressure:
	\begin{equation}
		\sum\limits_{j \in S_i}\D_{i,j}\Mpsi_{j}^{-1}\left(  \Q_{r(j),j} \bbpi_{r(j)}^{*} + \Q_{\ell(j),j} \bbpi_{\ell(j)}^{*}  \right) 
		= 
		\sum\limits_{j \in S_i}\D_{i,j}\Mpsi_{j}^{-1}\left(  \Gg_j 
		- \beta g \U_{r(j),j}\, \delta \Tm^* - \beta g \U_{\ell(j),j}\, \delta \Tm^*\right) 
		+\frac{1}{a_{ss} \Delta t} \sum\limits_{j \in S_i}\D_{i,j}\vv^* . \label{eq:pressure_IMEX} 
	\end{equation} 
	Once the values of $\vv^*,\bbpi^*$ and $\Tm^*$ are updated, the fluxes $\bar{\HH}^P_s$ and $\hat{\HH}^P_s$ can be evaluated as 
	\begin{eqnarray}
		\bar{\HH}^P_s=\left(
		\begin{array}{c}
			\Ldm(\hat{\HH}^{P,\vv}) \\
			\bar{\divop}(\alpha \hat{\gradop}(\bar{\temp}^*)) \\			
		\end{array}
		\right) \quad \mbox{ and } \quad
		\hat{\HH}^P_s=\left(
		\begin{array}{c}
			\Lmd(\hat{\divop}(\nu \bar{\gradop}(\vv^*)))- \hat{\gradop}(\bbpi^*) + \Mpsi_{j}^{-1}\mathbf{g}-\mathbf{g}\beta\Lmd(\delta \bar{\temp}^*) \\
			\Lmd(\bar{\HH}^{P,\temp})
		\end{array}
		\right).
	\end{eqnarray}
where $\bar{\divop}$ and $\hat{\gradop}$ refer to the discrete divergence and gradient operators which follows from \eqref{eq:ins.1}-\eqref{eq:ins.2} as
%\begin{gather}
%	\sum_{j \in S_i}\D_{i,j}\hbv_j^{n+1}=0 \label{eq:ins.1},\\
%	\Mpsi_{j} \TL_{h}^{\mathbf{v}}(\hbv_{j}) + \dt \, \left( \Q_{r(j),j} \bbpi_{r(j)}^{n+1}  + \Q_{\ell(j),j} \bbpi_{\ell(j)}^{n+1} \right) = 0, \label{eq:ins.2}
%\end{gather}
\begin{eqnarray}
	\bar{\divop}(\UUd^r)	\qquad &:& \qquad \bar{\divop}_i =	\sum_{j \in S_i}\D_{i,j}\UUd^r_j, \\
	\hat{\gradop}(\UUm^r) \qquad &:& \qquad  \hat{\gradop}_j= 	\Mpsi_{j}^{-1}\left(\Q_{r(j),j} \UUm^r_{r(j)}  + \Q_{\ell(j),j} \UUm^r_{\ell(j)}\right)	.
\end{eqnarray}

	After $\Ns$ stages, the final numerical solution at the new time level $t^{n+1}$ is updated using the coefficients $b_i$ of the Butcher tableaux. These are the same for both the explicit and the implicit tableau because we use the class of semi-implicit IMEX schemes proposed in \cite{BosFil2016}. Following \cite{TBP22}, we first update the trajectory
	\begin{equation}
		\xx^E_{j+1}=\xx^E_{j}-b_j \bar{\HH}^L_j(\xx) \qquad \forall j=1\ldots s,
	\end{equation}
	and the fluxes
	\begin{equation}
		\bar{\HH}^*_j	=\bar{\HH}^P_j(\xx+w_j\Delta t \bar{\HH}^L_j) \qquad \forall j=1\ldots s,
	\end{equation}	
	with $w_j   = \sum\limits_{r=1}^j a_{\Ns,r}-\sum\limits_{r=1}^{\Ns} b_r$. The final solution is then given by
	\begin{eqnarray}
		\bar{\UU}^{n+1}(\xx)=\bar{\UU}^n(\xx)+\sum\limits_{j=1}^{N_s}\Delta t \, b_j \, \bar{\HH}^*j 
	\end{eqnarray}
    on the main grid, and the same formulation is used for $\hat{\UU}^{n+1}(\xx)$ on the staggered mesh. More specifically, the previous equation is computed on the correct mesh where the variable is formally defined, i.e. either $\bar{\UU}^{n+1}$ or $\hat{\UU}^{n+1}$, and it is simply projected to the other mesh following \eqref{eq:statevecprj}.
The coupling between the semi-Lagrangian space-time discretization and the class of IMEX time integrators is not straightforward. Here, we have limited us to give an overview of the method in order to make the paper self-consistent, but the interested reader is referred to \cite{BT22} for all the related details and explanations on the construction of semi-Lagrangian IMEX schemes for both the non-conservative and the conservative case.
\subsection{Advection-diffusion}
As a simple particular example, we consider an advection-diffusion equation on a space-dependent velocity field. Also in this case, the transport part is solved explicitly using a semi-Lagrangian algorithm while the diffusion part is solved implicitly as discussed in Section \ref{sec:space_disc}. The advection-diffusion equation \eqref{eqn:advdiff} can then be written in the form \eqref{eq:cform} by setting
%\begin{equation}
%	\frac{\partial{C}}{\partial{t}}+\nabla \cdot \vv C= \nabla \cdot \left( \alpha \nabla C\right)
%\end{equation}
%that can be written in the previous general form by setting
$\UU=(C,\vv)$ and $\boldsymbol{\tau}=(1,0,0)^\top$ so that only the scalar quantity $C$ is evolved. Finally, $\mathbf{f}=(\vv C -\coeff \nabla C,0,0)$, $\mathbf{h}=0$ and $\mathbf{S}=0$. Due to the implicit discretization of the viscous terms, we have no pressure neither gravity contribution, thus the fluxes are simply given by
\begin{equation}
\bar{\HH}^P=\left(
\begin{array}{c}
	\bar{\divop}(\coeff \hat{\gradop}(C^*)) \\
	0
\end{array}
\right),
\end{equation}
that still requires the solution of a linear system for every IMEX stage and the evaluation of $C^*$ at the foot of the Lagrangian trajectory.
\section{Numerical experiments}\label{sec:num_exp}
We first test the novel high order semi-Lagrangian IMEX DG schemes for the advection-diffusion system. Then we will move to the incompressible Navier-Stokes equations and we will show applicability for density current computations in the case of small temperature gradients.

\subsection{Non-linear transport}
\begin{figure}[!htbp]
	\begin{center}
		\includegraphics[width=0.49\textwidth]{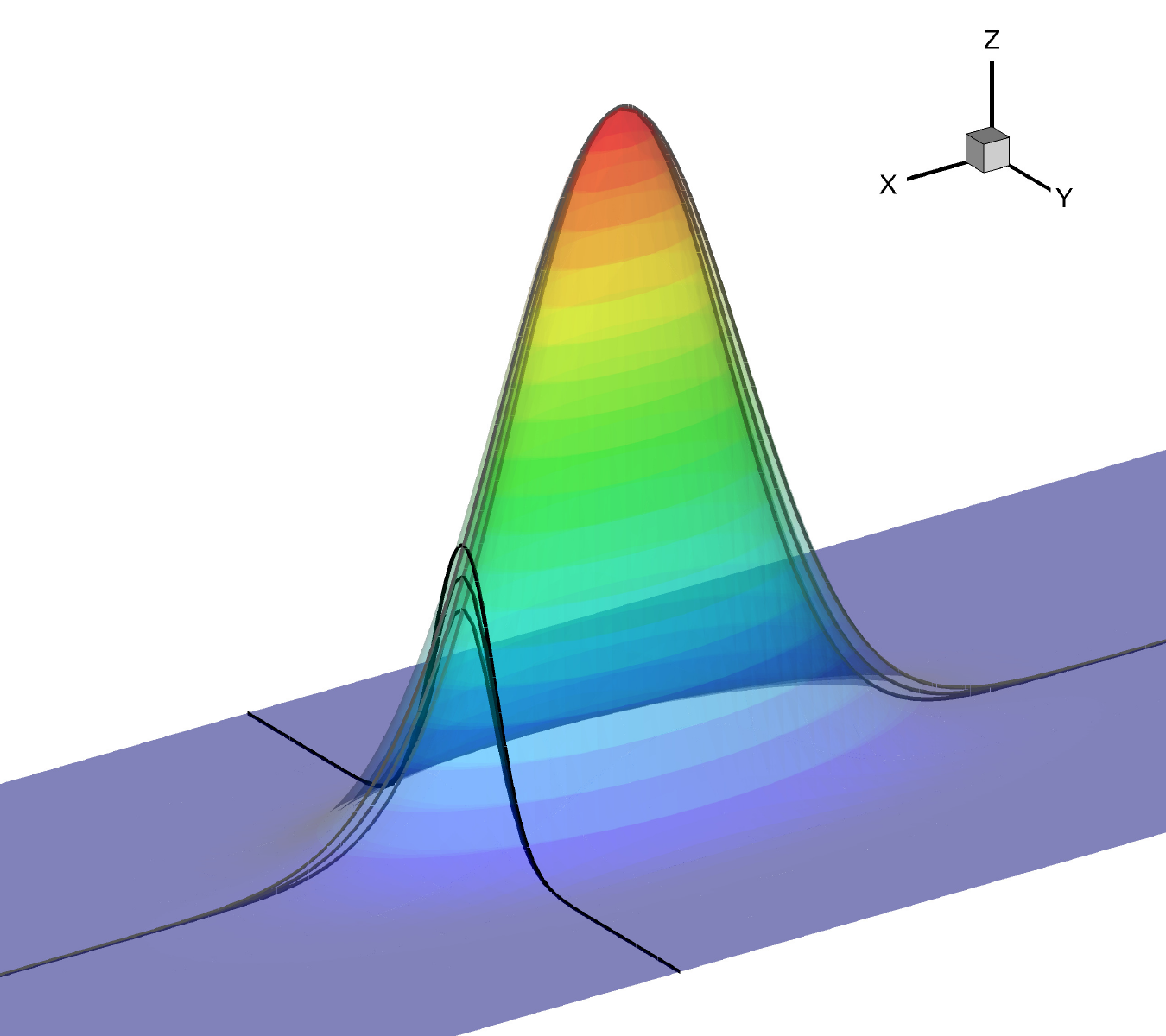}
		\caption{Two dimensional solution of non-linear transport test: comparison between exact and numerical solution with $R=1,2$ and extraction slices at $x=0.25$ and $y=0$.}
		\label{fig.NLT1}
	\end{center}
\end{figure}

\begin{figure}[!htbp]
	\begin{center}
		\includegraphics[width=0.49\textwidth]{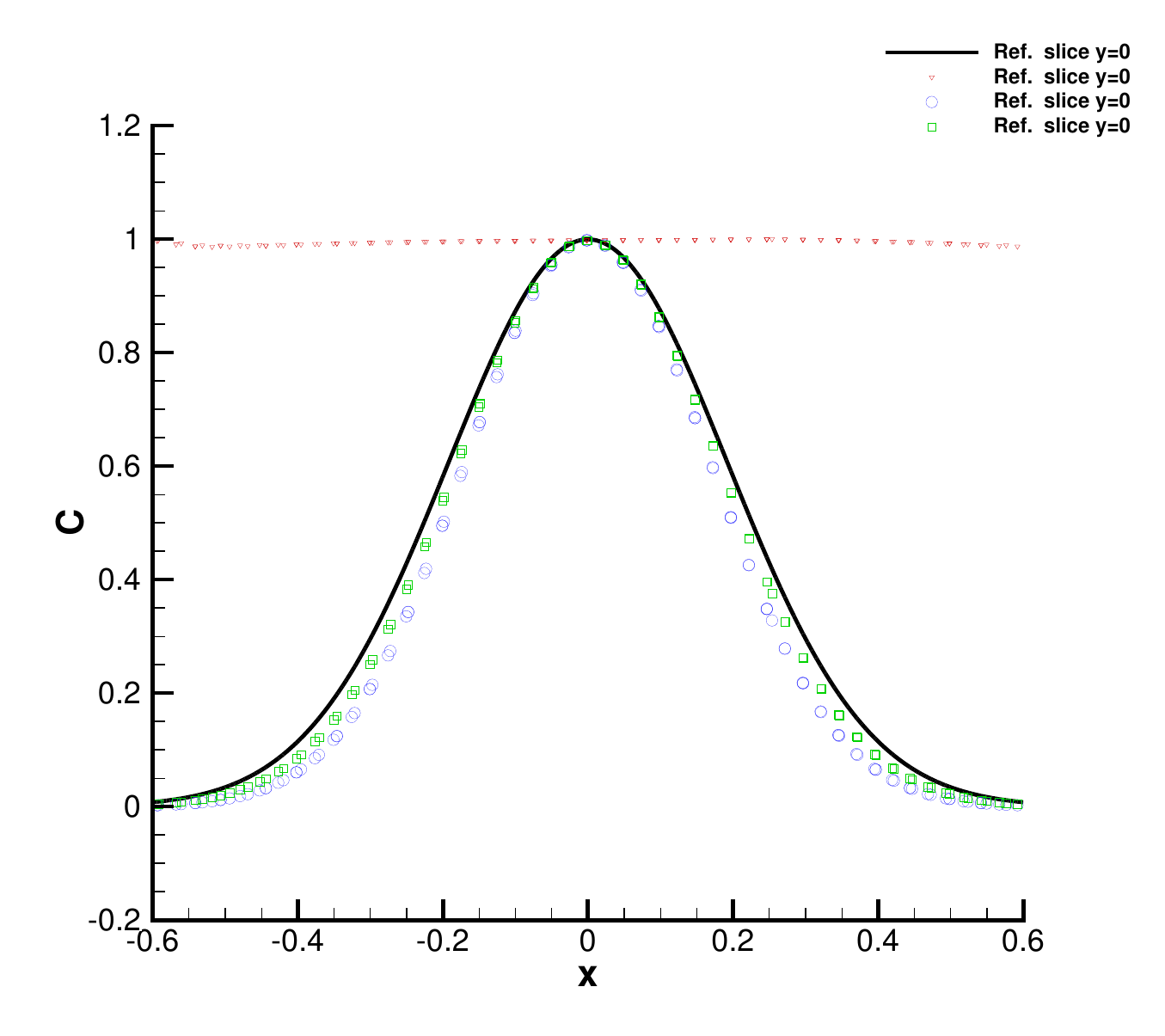}
		\includegraphics[width=0.49\textwidth]{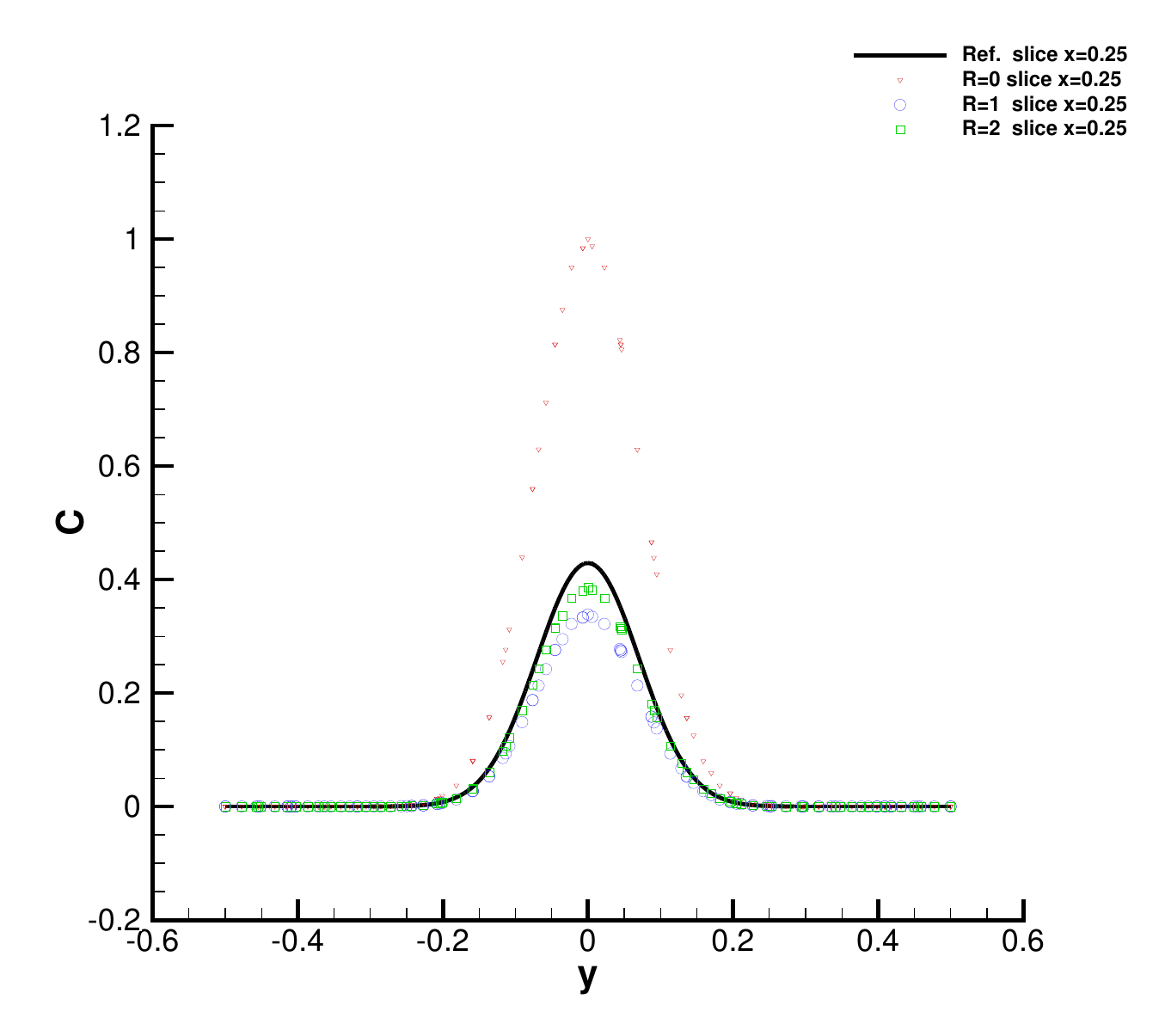}
		\caption{Non-linear transport. Comparison of the numerical solution obtained with $\Delta t=1.0$ and the exact one on the slice $y=0$ (left) and $x=0.25$ (right).}
		\label{fig.NLT2}
	\end{center}
\end{figure}

In this simple test we want to check the high order accuracy in time by solving the advection-diffusion equation with a non constant velocity field. For this test there is no evolution of the initial space-dependent velocity field $\vv(\xx,t)=\vv(\xx,0)=x$, and the pressure forces are neglected. The initial condition for the scalar quantity $C$ reads
\begin{equation}
	C(x,y,0)= e^{-100\sqrt{x^2+y^2}}.
\end{equation}
The time dependent exact solution can be found by following the characteristic equation, see \cite{TB19}:
\begin{equation}
	C(x,y,t)= exp \left(-100\sqrt{(xe^t)^2+y^2}\right).
\end{equation}
The diffusion coefficient is $\coeff=0$ and we set the final time to $t_{end}=1$. We consider a polynomial of degree $p=4$ on a domain $\Omega=[-2.5,2.5]\times [-0.5,0.5]$ covered with a main mesh composed by $\Ni=1124$ triangles. Periodic boundaries are imposed in the vertical direction and Dirichlet conditions in the $x$-direction. We start with a time step $\Delta t=1.0$ and then we successively refine both time and space. A visual comparison between the exact solution and the one computed with the IMEX schemes $R=0,1,2$ and a single step $\Delta t =1.0$ is reported in Figure $\ref{fig.NLT1}$, while Figure $\ref{fig.NLT2}$ shows a comparison in terms of one-dimensional cuts. The resulting convergence is reported in Table $\ref{tab:ct1}$, confirming that the formal order of convergence is achieved for all the adopted IMEX schemes (see \ref{app.imex}). 

\begin{table}[!h]
	\begin{center}
		\begin{tabular}{|c|c||c|c|c|c|c|c|}
			\hline
			\multicolumn{2}{|c||}{Ref.}    & \multicolumn{2}{c|}{$R=0$} & \multicolumn{2}{c|}{$R=1$} & \multicolumn{2}{c|}{$R=2$}
			\\\hhline{|--|-|-|-|-|-|-|}
			$\Delta t$& $r_{in}$ & $L_2$-error & order  & $L_2$-error & order   &  $L_2$-error & order  \\\hline\hline
			$1.0$  		&$2.14E-02$&  $7.406666034E-01$     & $-$    & $2.157102724E-02$ & $-$   & $1.045342032E-02$ & $-$  \\\hline
			$0.5$   	&$1.07E-02$&  $7.511722708E-02$     & $3.3$  & $4.689157993E-03$ & $2.2$ & $7.068162565E-04$ & $3.9$  \\\hline
			$0.33$   	&$7.14E-03$&  $4.070014613E-02$     & $1.5$  & $1.942549475E-03$ & $2.2$ & $1.664730210E-04$ & $3.6$  \\\hline
			$0.25$     &$5.36E-03$	&  $2.794974880E-02$     & $1.3$  & $1.052023386E-03$ & $2.1$ & $6.247319230E-05$ & $3.4$  \\\hline
		\end{tabular}
	\end{center}	
	\caption{Convergence results for the non-linear transport problem.}\label{tab:ct1}
\end{table}	

\subsection{Advection-diffusion}
As explained in Section \ref{sec:SLIMEX}, the semi-Lagrangian IMEX-DG schemes need to properly transport the fluxes according to the velocity field in order to be compatible with the stages of the pure IMEX scheme. In order to check this consistency we perform a simple but not trivial test involving both advection and diffusion. We first set a constant velocity field $\vv(\xx,t)=(0.2,0.0)$ and a diffusion coefficient $\coeff=10^{-3}$ on a domain $\Omega=[-2.5,2.5]\times[-0.5,0.5]$ initially covered with $\Ni=1124$ triangles. The exact solution for the evolution of a step function over a constant velocity field is the simple Lagrangian transport of the solution of the stationary heat equation and reads:
\begin{equation}
	C(x,y,t)=\frac{1}{2}-\frac{1}{2} erf \left( \frac{x-x_0-ut}{\sqrt{4 \coeff (t +t_0)}} \right),
\end{equation}
with $x_0=-0.5$ and $t_0=0.5$ in order to avoid oscillations due to the $L_2$ projection of the initial condition at $t=0$.
We take $p=4$ and we sequentially reduce both the space and the time mesh size. We run several IMEX schemes $R=0,1,2$ up to $t_{end}=2.0s$. As for the previous case, we consider periodic boundaries in $y$-direction and Dirichlet conditions in $x$-direction where the exact solution is prescribed.
\begin{figure}[!htbp]
	\begin{center}
		\includegraphics[width=0.6\textwidth]{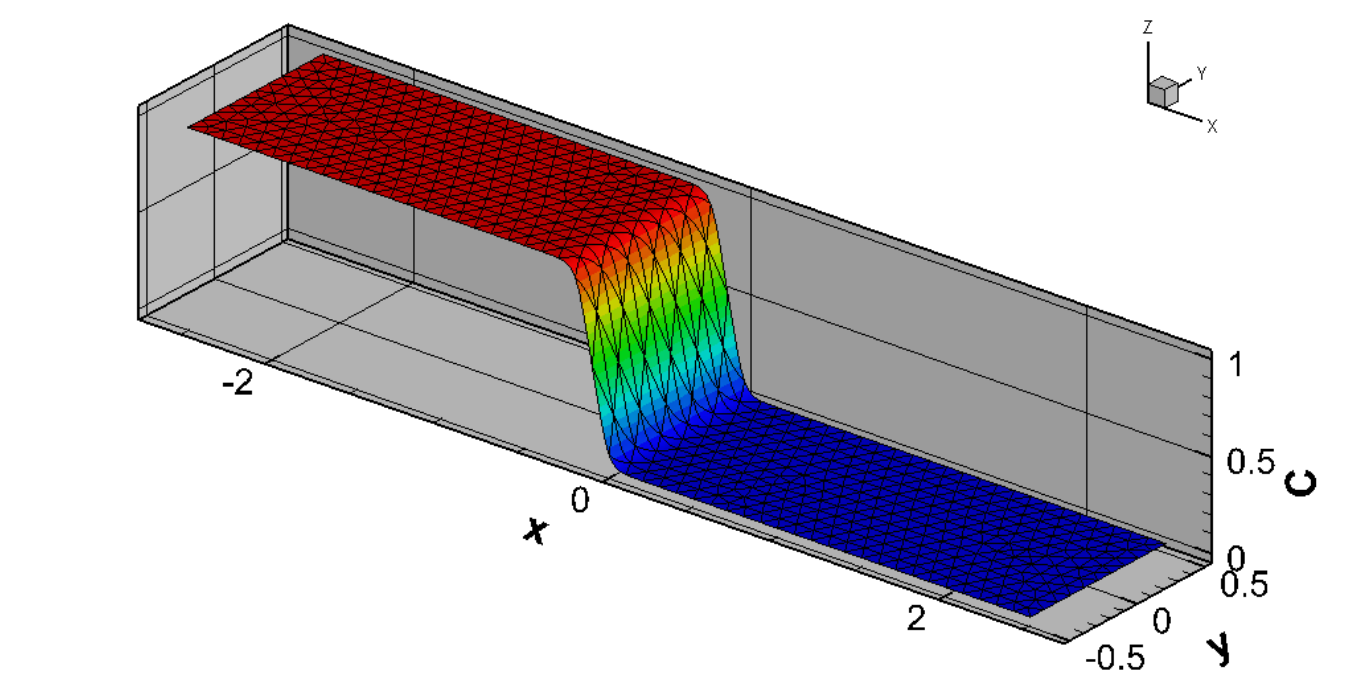}
		\caption{Advection-diffusion test. Plot of the solution at time $t_{end}=2$ with $R=2$ and the initial mesh with $\Ni=1124$.}
		\label{fig.AD1}
	\end{center}
\end{figure}
A sketch of the solution obtained using $R=2$ at the first refinement is shown in Figure \ref{fig.AD1}.
\begin{table}[!htbp]
	\begin{center}
		\begin{tabular}{|c|c||c|c|c|c|c|c|}
			\hline
			\multicolumn{2}{|c||}{Ref.}    & \multicolumn{2}{c|}{$R=0$} & \multicolumn{2}{c|}{$R=1$} & \multicolumn{2}{c|}{$R=2$}
			\\\hhline{|--|-|-|-|-|-|-|}
			$\Delta t$& $r_{in}$ & $L_2$-error & order  & $L_2$-error & order   &  $L_2$-error & order  \\\hline\hline
			$2.00$  		&$2.14E-02$&  $1.175200075E-02$     & $-$    & $4.155363690E-03$ & $-$   & $2.246920883E-03$ & $-$  \\\hline
			$1.00$   	&$1.07E-02$&  $6.685962372E-03$     & $0.8$  & $5.856947865E-04$ & $2.8$ & $2.317034705E-04$ & $3.3$  \\\hline
			$0.67$   	&$7.14E-03$&  $4.671479996E-03$     & $0.9$  & $2.449695824E-04$ & $2.1$ & $7.614394839E-05$ & $2.7$  \\\hline
			$0.50$      &$5.36E-03$&  $3.589321174E-03$     & $0.9$  & $1.349022601E-04$ & $2.1$ & $3.450490896E-05$ & $2.8$  \\\hline
		\end{tabular}
	\end{center}	
	\caption{Convergence results for the advection-diffusion problem.}\label{tab:ct2}
\end{table}	
The resulting convergence rates are reported in Table $\ref{tab:ct2}$ for all the cases $R=0,1,2$ and time partitions $\Delta t = \{t_{end}/1,t_{end}/2,t_{end}/3,t_{end}/4\}$. The expected order of accuracy can be observed also in this case for all the adopted IMEX schemes. To be noted that the polynomial degree is $4$, hence ensuring that the convergence rate is limited only by the time accuracy of the method. In all the cases the diffusion system can be solved using the conjugate gradient method. This is not the case when a fully space-time is adopted, see \cite{TD14,TD15,TD16}, where the high order discretization in time fails to maintain the system symmetric, even at second order.

\subsection{Advection-diffusion with non-constant velocity field}
\begin{figure}[!htbp]
	\begin{center}
		\includegraphics[width=0.9\textwidth]{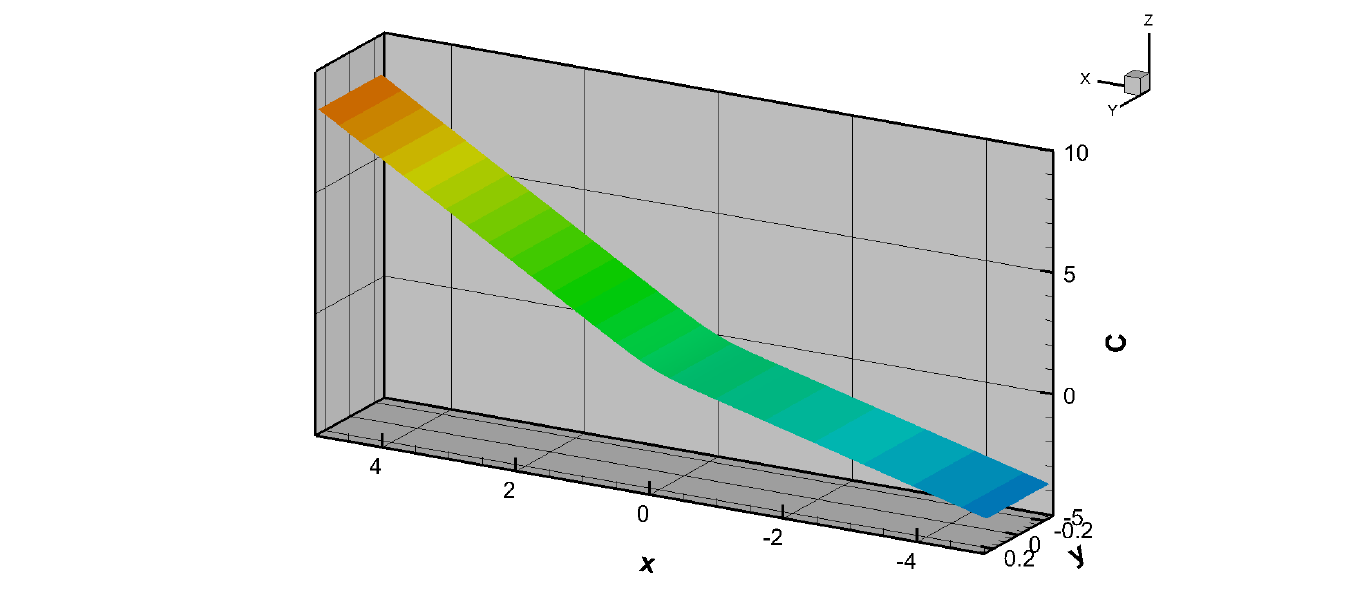}
		\caption{Advection-diffusion with non-constant velocity field. Solution of the scalar quantity $C$ at time $t_{end}=0.2$.}
		\label{fig.nlad1}
	\end{center}
\end{figure}

The previous test involved the combination of a constant velocity field and viscous effects. 
We test here our numerical algorithm in a case where both the viscous and non-constant transport are considered. As derived in \cite{TBP22}, an exact solution of the scalar equation 
\begin{equation}
	\diff{C}{t}-x\diff{C}{x} = \coeff \frac{\partial^2 C}{\partial x^2},
	\label{AD2.1}
\end{equation} 
reads
\begin{equation}
	C_{ex}(x,t)=e^{t}\left[ -x-\frac{\sqrt{2\pi\coeff}}{2} erf(\frac{1}{\sqrt{2}x})x-e^{-\frac{x^2}{2\coeff}}\coeff \right].
\end{equation}
The computational domain is $\Omega=[-5,5] \times [-0.25,0.25]$ which is then extended up to $10$ in the $x$-direction in order to avoid boundary interference. The velocity field reads $\vv=(-x,0)$ according to \eqref{AD2.1}, and we choose the viscosity parameter $\coeff=0.1$ and the final time $t_{end}=0.2$. As for the previous tests, we sequentially refine the mesh starting from a triangulation that counts $\Ni=1536$  cells and using only one time step, i.e. $\Delta t=t_{end}$. We use $p=4$ and the exact solution as initial condition. A sketch of the solution obtained at the final time is depicted in Figure \ref{fig.nlad1} while the resulting $L_2$ error on $\Omega$ is reported in Table \ref{tab:ct3}. Even for this advection-diffusion test with space-dependent velocity field the optimal order of convergence is reached.
\begin{table}[!h]
	\begin{center}
		\begin{tabular}{|c|c||c|c|c|c|c|c|}
			\hline
			\multicolumn{2}{|c||}{Ref.}    & \multicolumn{2}{c|}{$R=0$} & \multicolumn{2}{c|}{$R=1$} & \multicolumn{2}{c|}{$R=2$}
			\\\hhline{|--|-|-|-|-|-|-|}
			$\Delta t$& $r_{in}$ & $L_2$-error & order  & $L_2$-error & order   &  $L_2$-error & order  \\\hline\hline
			$0.2$  		&$2.80E-02$&  $1.485813636E-01$     & $-$    & $5.111092318E-03$ & $-$   & $1.508897625E-04$ & $-$  \\\hline
			$0.1$   	&$1.40E-02$&  $7.916125589E-02$     & $0.9$  & $1.344222242E-03$ & $1.9$ & $2.199110494E-05$ & $2.8$  \\\hline
			$0.067$   	&$9.35E-03$&  $5.396364823E-02$     & $0.9$  & $6.075716380E-04$ & $2.0$ & $6.686687647E-06$ & $2.9$  \\\hline
			$0.05$      &$7.01E-03$&  $4.093565054E-02$     & $1.0$  & $3.446405417E-04$ & $2.0$ & $2.859700682E-06$ & $3.0$  \\\hline
		\end{tabular}
	\end{center}	
	\caption{Convergence results for the advection-diffusion test with non-constant velocity field.}\label{tab:ct3}
\end{table}	

\subsection{Taylor Green vortex}
\begin{figure}[!htbp]
	\begin{center}
		\includegraphics[width=0.4\textwidth]{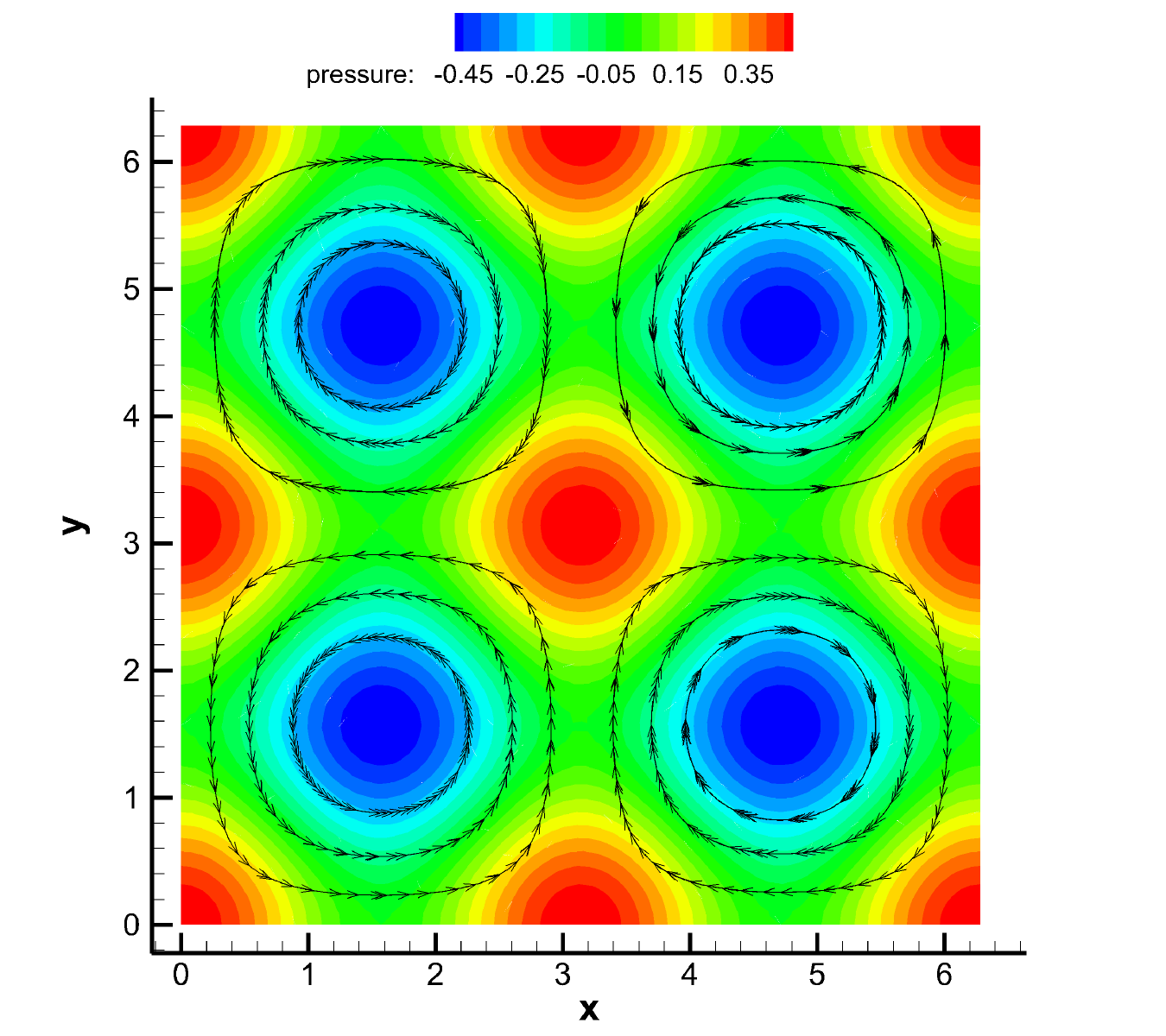}
		\includegraphics[width=0.4\textwidth]{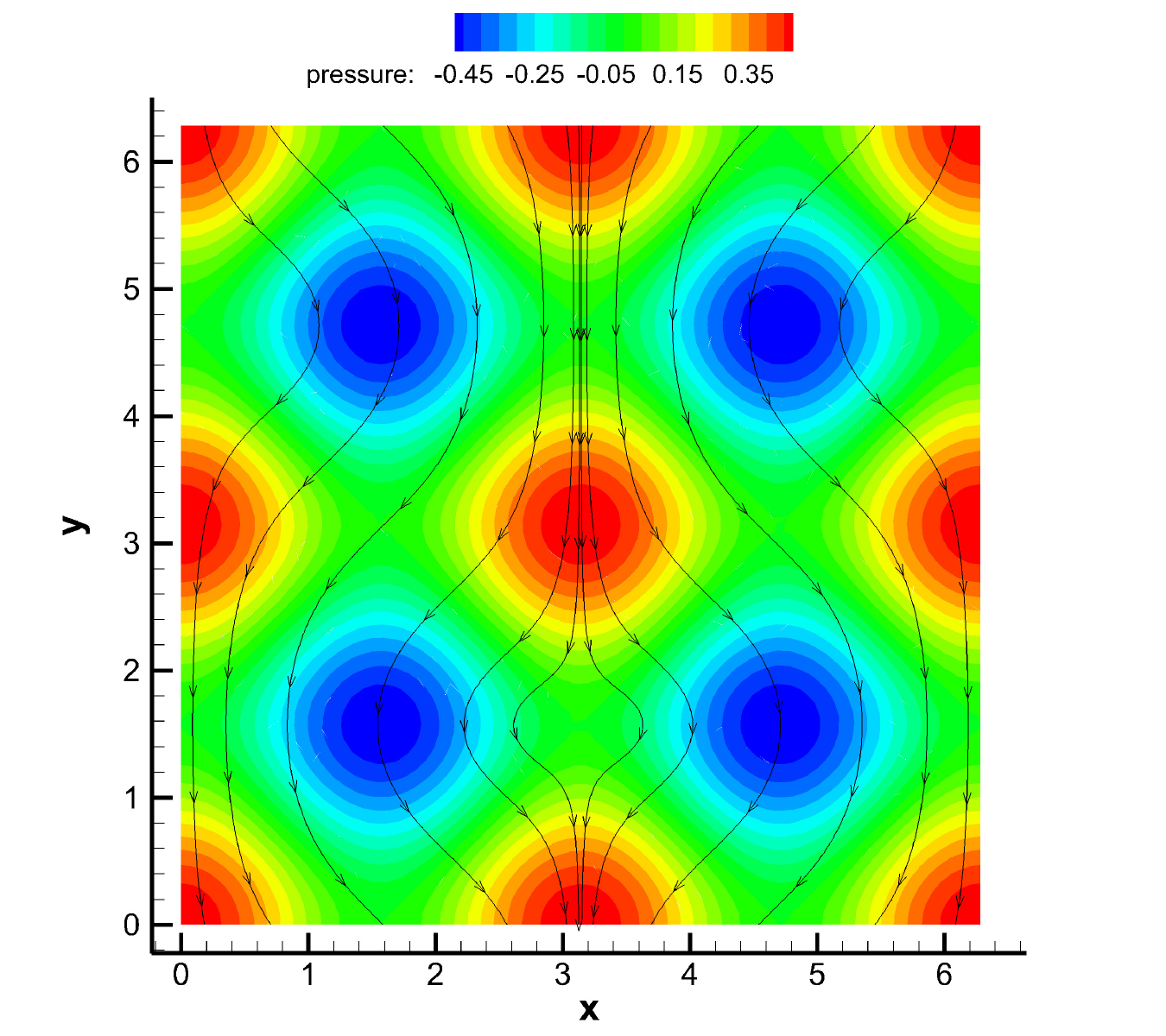}
		\caption{Taylor-Green vortex. Pressure and velocity field at time $t=0$ (left) and $t=t_{end}$ (right).}
		\label{fig.TGV.1}
	\end{center}
\end{figure}

Here, we consider the two-dimensional incompressible Navier-Stokes model. We numerically solve the two-dimensional Taylor Green vortex problem that is a benchmark for incompressible fluids. The exact solution of this test is given by the exact balance between kinetic and pressure forces and reads:
\begin{eqnarray}
	u_{ex}(x,y,t)&=& \sin(x)\cos(y)e^{-2\nu t}, \label{eq:tgvsol1}\\
	v_{ex}(x,y,t)&=& -\cos(x)\sin(y)e^{-2\nu t} +gt, \\
	p_{ex}(x,y,t)&=& \frac{1}{4}(\cos(2x)+\cos(2y))e^{-2\nu t}. \label{eq:tgvsol3}
\end{eqnarray}
on a domain $\Omega=[0,2\pi]^2$. By considering a non zero gravity force $g=-9.81$, the definition \eqref{eq:tgvsol1}-\eqref{eq:tgvsol3} remains a solution only if we add an additional source contribution 
\begin{equation}
	S=\left( 
	\begin{array}{c}
		0 \\ -gt\sin(x)\sin(y)e^{-2\nu t} \\ gt \cos(x)\cos(y)e^{-2\nu t} \\ 0
	\end{array}
	\right)
\end{equation}
in the momentum equation. The initial mesh consists in only $\Ni=116$ triangles, we set $p=R$ and $\nu=10^{-3}$. The fluid is initially assigned the exact solution with $t=0$. The final time is taken as $t_{end}=0.1$ and we carry out the simulation with several time steps in order to show time convergence. Periodic boundaries are considered in all the directions and, since the velocity field is accelerating in the $y$ directions, the trajectories are expected to pass through the vertical periodic boundary.
The pressure and velocity field are depicted in Figure $\ref{fig.TGV.1}$ at the initial time $t=0$ and at the final time $t=t_{end}$ for the case $R=2$ and $\dt=0.025$. As expected, the trajectories that are initially closed start to pass through the periodic boundary according to the increasing vertical velocity.
The resulting convergence studies are reported in Table $\ref{tab:tgv1}$, where the $L_2$ errors for the velocity vector field are defined as
\begin{equation}
	L_2^{err}=\sqrt{\int\limits_{\Omega}(u-u_{ex})^2+(v-v_{ex})^2}.
\end{equation}

% Old table using N_i=902 elements
%\begin{table}[H]
%	\begin{center}
%		\begin{tabular}{|c|c||c|c|c|c|c|c|}
%			\hline
%			\multicolumn{2}{|c||}{Ref.}    & \multicolumn{2}{c|}{$R=0$} & \multicolumn{2}{c|}{$R=1$} & \multicolumn{2}{c|}{$R=2$}
%			\\\hhline{|--|-|-|-|-|-|-|}
%			$\Delta t$& $r_{in}$ & $L_2$-error & order  & $L_2$-error & order   &  $L_2$-error & order  \\\hline\hline
%			$0.1$  		&$2.14E-02$&  $3.128594177E-02$     & $-$    & $2.933260274E-02$ & $-$   & $1.748674971E-03$ & $-$  \\\hline
%			$0.05$   	&$1.07E-02$&  $1.067774760E-02$     & $1.6$  & $7.942026218E-03$ & $1.9$ & $2.385796939E-04$ & $2.9$  \\\hline
%			$0.025$   	&$7.14E-03$&  $6.507963115E-03$     & $1.2$  & $3.598252856E-03$ & $2.0$ & $7.250558017E-05$ & $2.9$  \\\hline
%			$0.0125$     &$5.36E-03$	& $4.797432952E-03$     & $1.1$  & $2.061560492E-03$ & $1.9$ & $3.076933821E-05$ & $3.0$  \\\hline
%		\end{tabular}
%	\end{center}	
%	\caption{Convergence rate obtained for the Taylor-Green vortex with gravity and source terms.}\label{tab:tgv1}
%\end{table}	
\begin{table}[!h]
	\begin{center}
		\begin{tabular}{|c|c||c|c|c|c|c|c|}
			\hline
			\multicolumn{2}{|c||}{Ref.}    & \multicolumn{2}{c|}{$R=0$} & \multicolumn{2}{c|}{$R=1$} & \multicolumn{2}{c|}{$R=2$}
			\\\hhline{|--|-|-|-|-|-|-|}
			$\Delta t$& $r_{in}$ & $L_2$-error & order  & $L_2$-error & order   &  $L_2$-error & order  \\\hline\hline
			$0.1$  		&$2.51E-01$		&  $1.206001315E+00$     & $-$    & $1.960961356E-01$ & $-$   & $2.208824449E-02$ & $-$  \\\hline
			$0.05$   	&$1.25E-01$		&  $6.188176625E-01$     & $1.0$  & $5.778024572E-02$ & $1.8$ & $3.397912274E-03$ & $2.7$  \\\hline
			$0.033$   	&$8.39E-02$		&  $4.319855848E-01$     & $0.9$  & $2.833213299E-02$ & $1.8$ & $1.136618866E-03$ & $2.7$  \\\hline
			$0.025$    &$6.29E-02$	    &  $3.463894119E-01$     & $0.8$  & $1.719363340E-02$ & $1.7$ & $5.226162393E-04$ & $2.7$  \\\hline
		\end{tabular}
	\end{center}	
	\caption{Convergence rates obtained for the Taylor-Green vortex with gravity and source terms.}\label{tab:tgv1}
\end{table}	
The expected order of convergence can be observed also for this test case. We underline that, since the DG discretization is adopted in space only and high order in time is achieved using IMEX schemes, the resulting linear system for the pressure is symmetric and it can thus be solved using the conjugate gradient method. On the contrary, genuinely high order space-time DG schemes would require a more general and slower GMRES solver due to the non-symmetry of the space-time operators. 

To analyze the performance of the new SL IMEX-DG schemes at different viscosity scales, let us now consider different values of the viscosity coefficient that runs from $\nu=10^{-6}$ up to $\nu=10^{-2}$ and we set $g=0$ so that no artificial source terms are needed. The resulting convergence rates for $R=0,1,2$ are reported in Table \ref{tab:tgv2}. An order of convergence of $p+1/2$ is observed for small values of $\nu$ that is in accordance to what observed in \cite{NatConv20} using a fully space time DG method. The remaining convergence rates are in accordance with the formal order of accuracy of the method.
\begin{table}[!h]
	\begin{center}
		\begin{tabular}{|c||c|c|c|c|c|c|}
			\hline
			$R=2$ & \multicolumn{2}{c|}{$\nu=1e-6$} & \multicolumn{2}{c|}{$\nu=1e-4$} & \multicolumn{2}{c|}{$\nu=1e-2$}
			\\\hhline{|-|-|-|-|-|-|-|}
			Ref. & $L_2$-error & order  & $L_2$-error & order   &  $L_2$-error & order  \\\hline\hline
			$1$	&  $2.209894331E-02$     & $-$    & $2.209780979E-02$ & $-$   & $2.203577143E-02$ & $-$  \\\hline
			$2$ &  $3.531993332E-03$     & $2.6$  & $3.528710361E-03$ & $2.6$ & $3.336958270E-03$ & $2.7$  \\\hline
			$3$ &  $1.277805244E-03$     & $2.5$  & $1.273844779E-03$ & $2.5$ & $1.083373195E-03$ & $2.8$  \\\hline
			$4$ &  $6.575808369E-04$     & $2.3$  & $6.529412708E-04$ & $2.3$ & $4.739906301E-04$ & $2.9$  \\\hline
			\hline
			$R=1$ & \multicolumn{2}{c|}{$\nu=1e-6$} & \multicolumn{2}{c|}{$\nu=1e-4$} & \multicolumn{2}{c|}{$\nu=1e-2$}
			\\\hhline{|-|-|-|-|-|-|-|}
			Ref. & $L_2$-error & order  & $L_2$-error & order   &  $L_2$-error & order  \\\hline\hline
			$1$	&  $1.961775016E-01$     & $-$    & $1.961693618E-01$ & $-$   & $1.954327805E-01$ & $-$  \\\hline
			$2$ &  $5.782369642E-02$     & $1.8$  & $5.779291996E-02$ & $1.8$ & $5.526283321E-02$ & $1.8$  \\\hline
			$3$ &  $2.872724184E-02$     & $1.7$  & $2.868131625E-02$ & $1.7$ & $2.553257902E-02$ & $1.9$  \\\hline
			$4$ &  $1.763853734E-02$     & $1.7$  & $1.758390564E-02$ & $1.7$ & $1.450076476E-02$ & $2.0$  \\\hline
			\hline
			$R=0$ & \multicolumn{2}{c|}{$\nu=1e-6$} & \multicolumn{2}{c|}{$\nu=1e-4$} & \multicolumn{2}{c|}{$\nu=1e-2$}
			\\\hhline{|-|-|-|-|-|-|-|}
			Ref. & $L_2$-error & order  & $L_2$-error & order   &  $L_2$-error & order  \\\hline\hline
			$1$	&  $1.206014673E+00$     & $-$    & $1.206013343E+00$ & $-$   & $1.205886729E+00$ & $-$  \\\hline
			$2$ &  $6.094749600E-01$     & $1.0$  & $6.094875276E-01$ & $1.0$ & $6.107868307E-01$ & $1.0$  \\\hline
			$3$ &  $4.073098796E-01$     & $1.0$  & $4.073271692E-01$ & $1.0$ & $4.091297136E-01$ & $1.0$  \\\hline
			$4$ &  $3.059365751E-01$     & $1.0$  & $3.059564052E-01$ & $1.0$ & $3.080442722E-01$ & $1.0$  \\\hline
		\end{tabular}
	\end{center}	
	\caption{Convergence rates obtained for the Taylor-Green vortex with different viscosity coefficients.}\label{tab:tgv2}
\end{table}	
\subsection{Rising warm bubble}
\begin{figure}[!htbp]
	\begin{center}
		\includegraphics[width=0.30\textwidth]{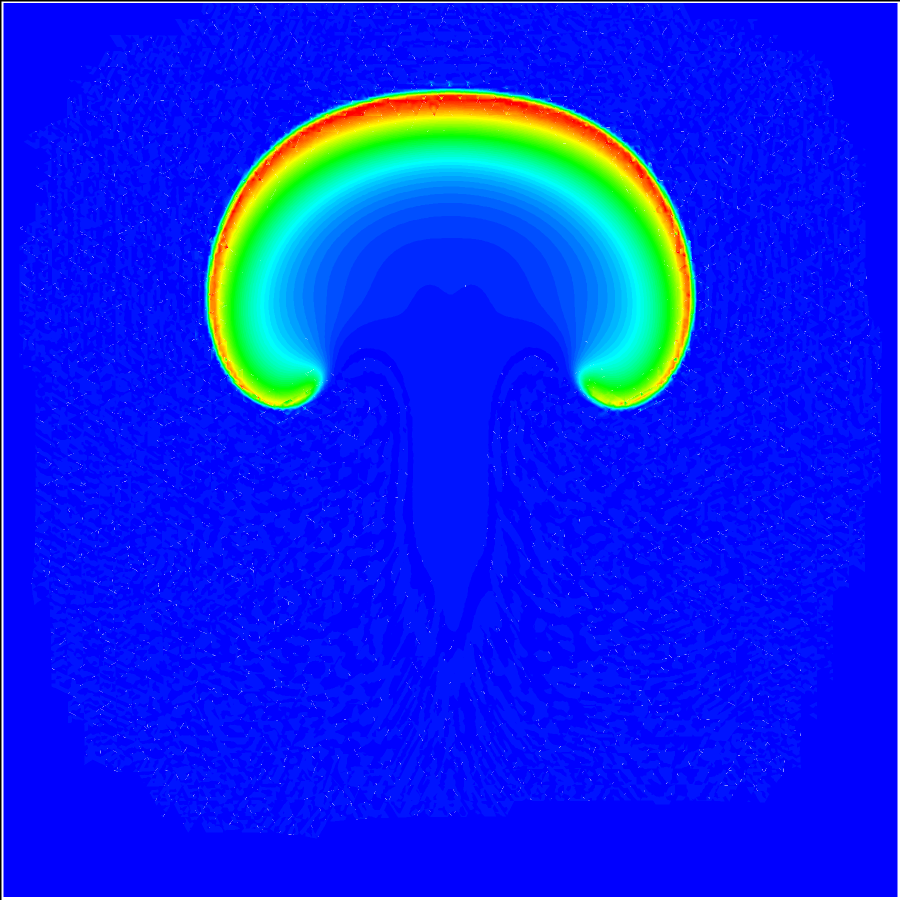}
		\includegraphics[width=0.30\textwidth]{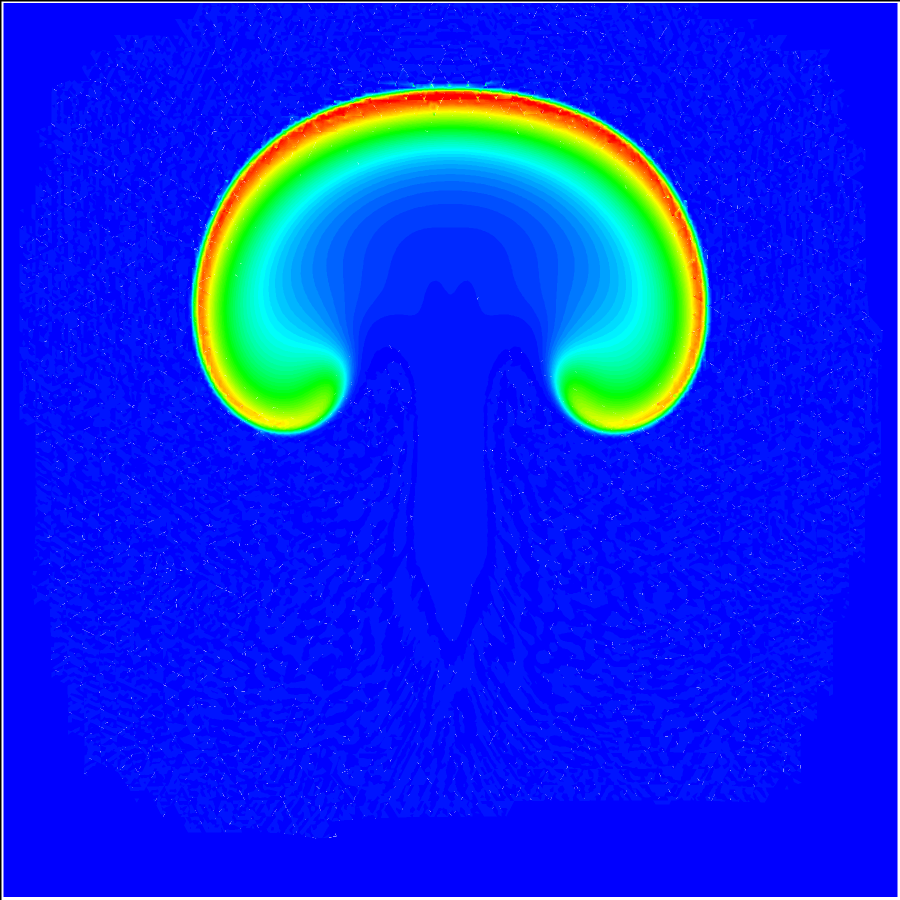}
		\includegraphics[width=0.30\textwidth]{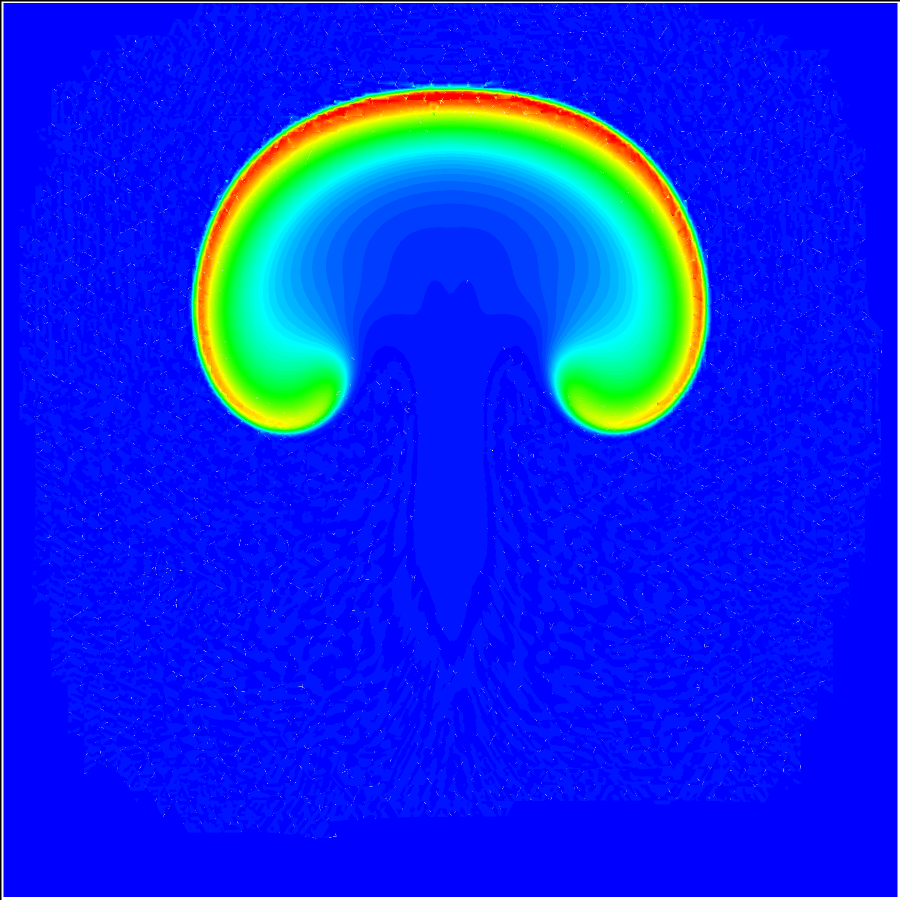} \\
		\includegraphics[width=0.30\textwidth]{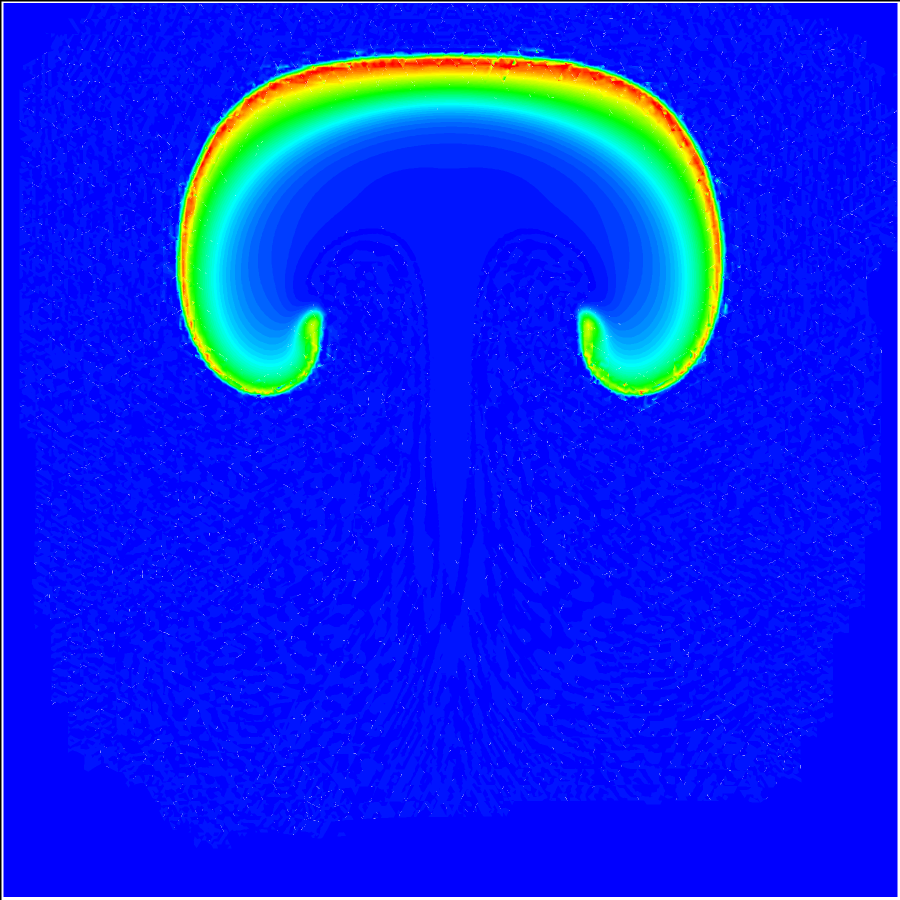}
		\includegraphics[width=0.30\textwidth]{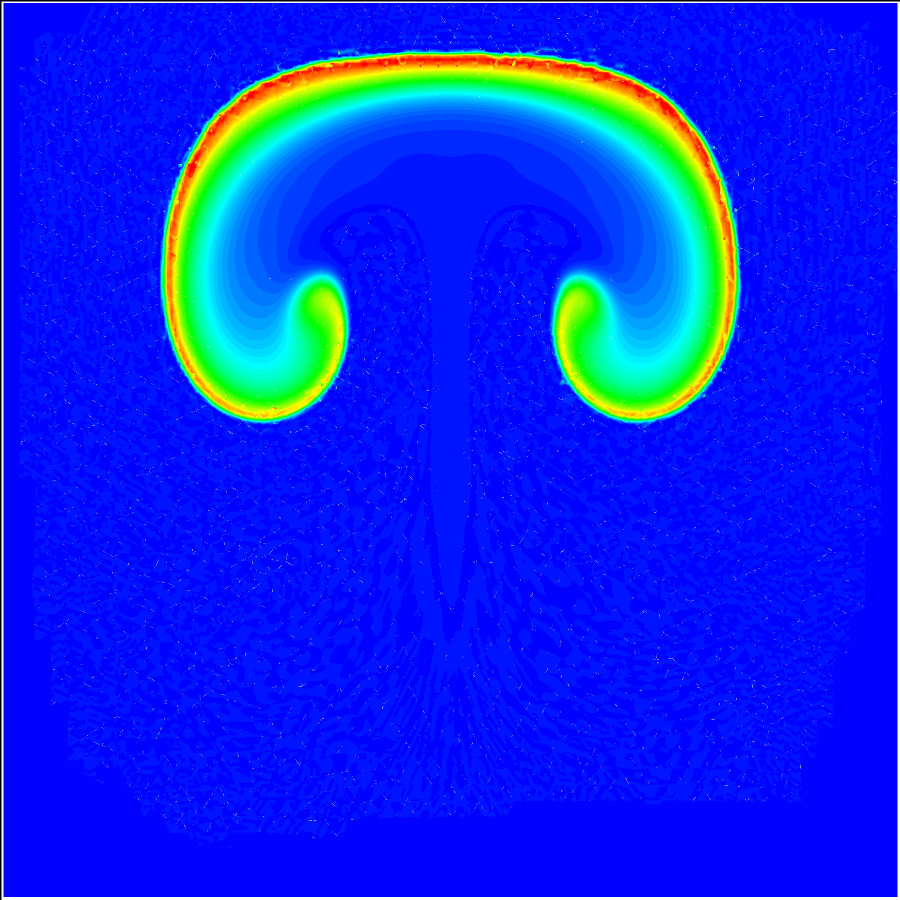}
		\includegraphics[width=0.30\textwidth]{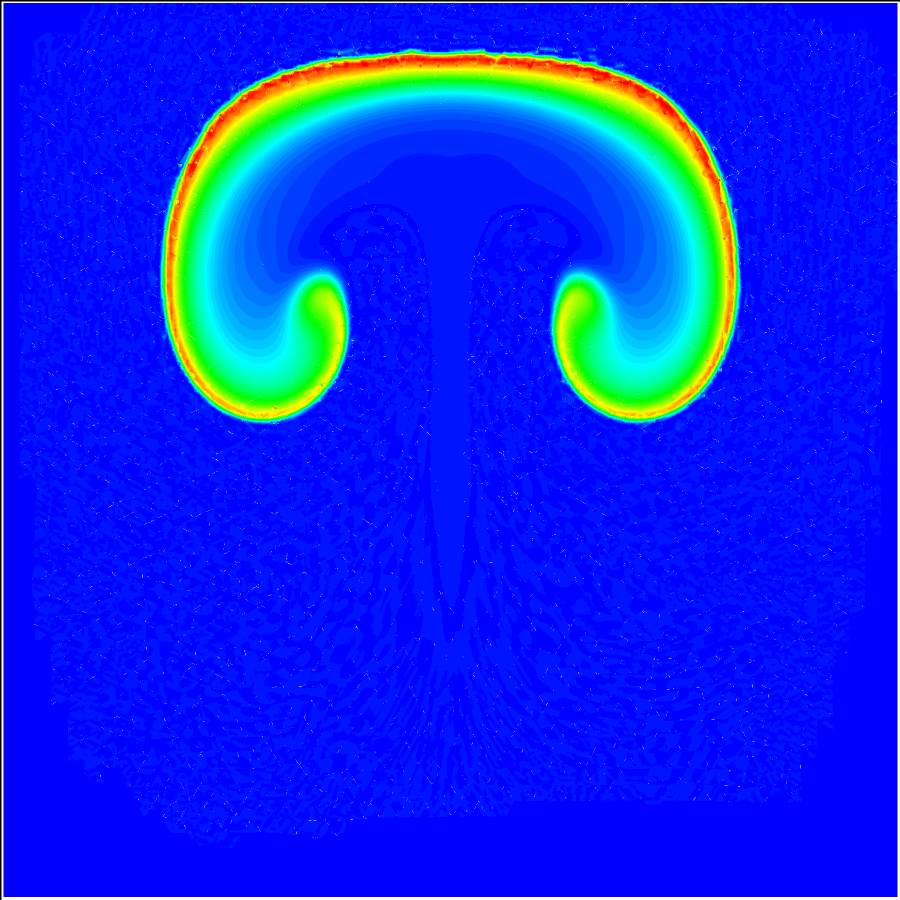} \\
		\includegraphics[width=0.30\textwidth]{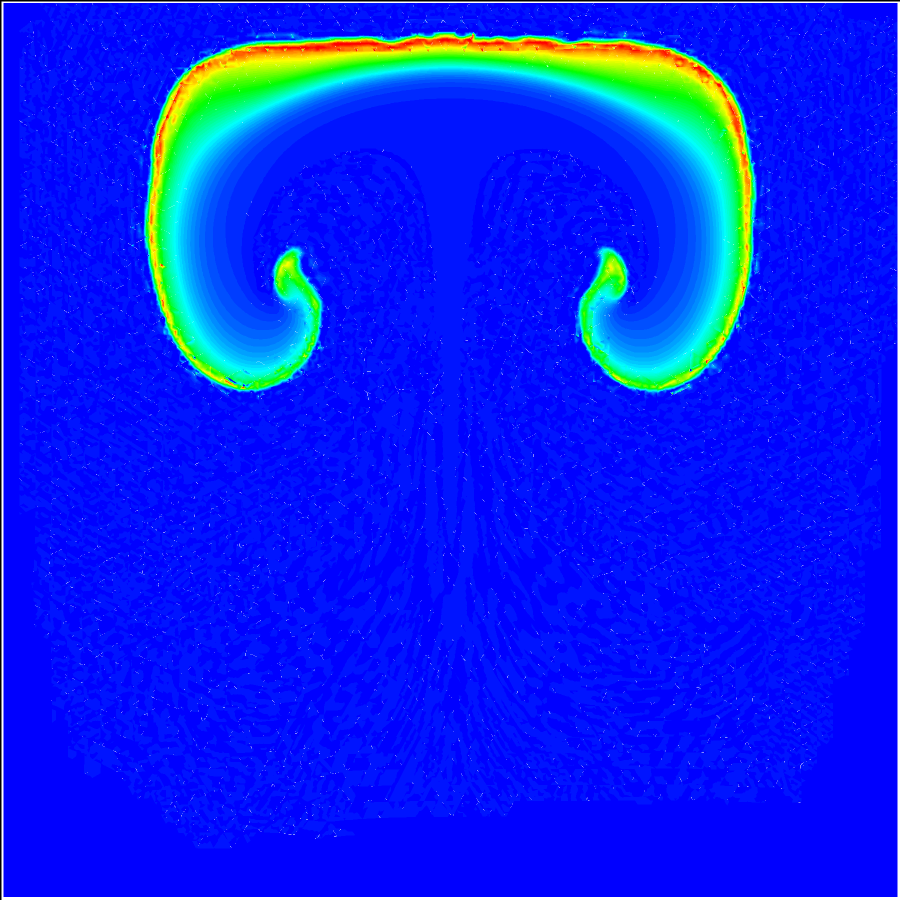}
		\includegraphics[width=0.30\textwidth]{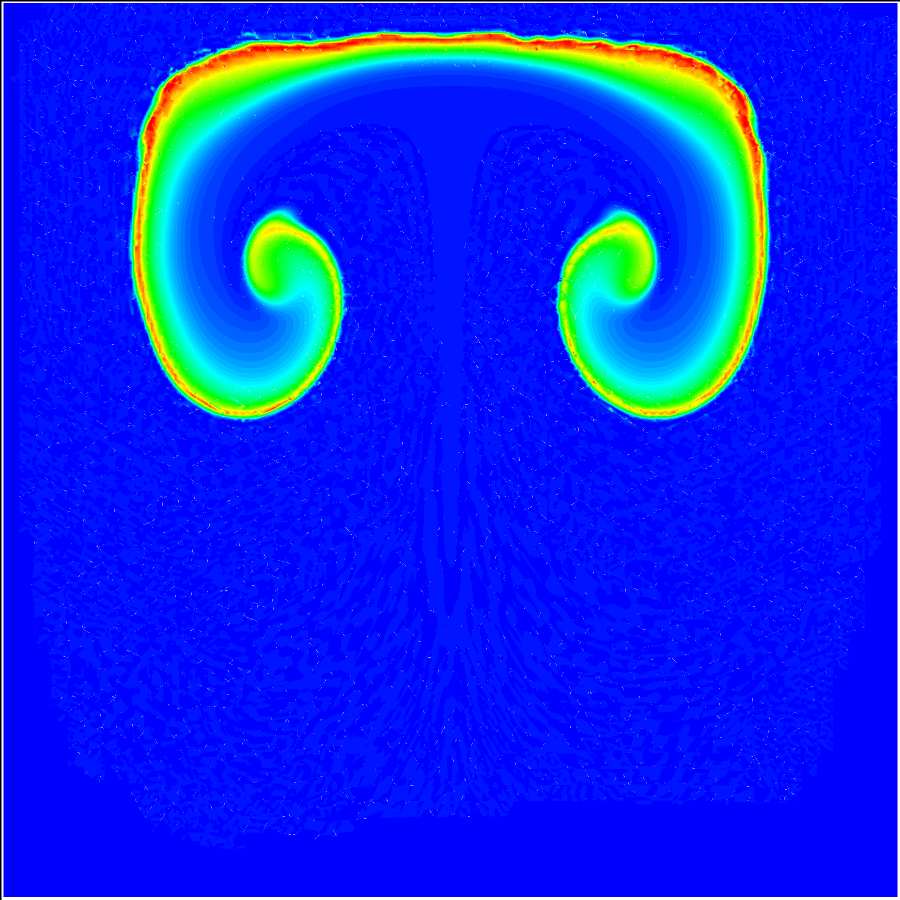}
		\includegraphics[width=0.30\textwidth]{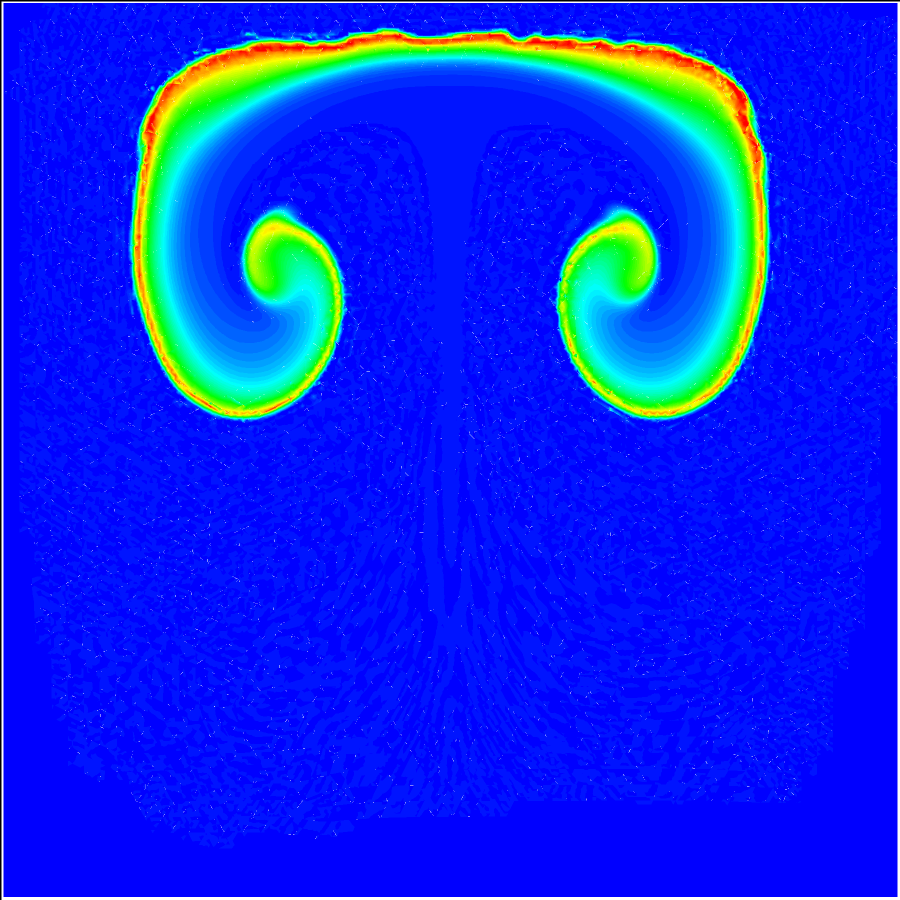}
		\caption{Temperature profile for the warm bubble test at times $t=600,700,800$ from the first to the last row and with $R=0,1,2$ from the first to the last column, respectively.}
		\label{fig.WV.1}
	\end{center}
\end{figure}

This test cases involves both the Navier-Stokes model and the advection-diffusion evolution of a scalar quantity. Under the Boussinesq approximation, the coupling can simply be done using a source term in the momentum equation and hence our SL IMEX-DG algorithm can be applied by setting $\UU=(p,\vv,\temp)$ and $\gamma=(0,1,1)$ according to equation \eqref{eq:imex_NC}. The initial condition for the quantity $\temp$ is taken from \cite{NatConv20} and it reads
\begin{eqnarray}
	\temp(x,0)=\left\{
	\begin{array}{lc}
		0 & r>r_b \\
		\frac{\temp_b}{2}\left[1 + \cos\left(\frac{\pi r}{\temp_b}\right)\right] & r \leq r_b
	\end{array}\right. ,
\end{eqnarray}
where $(x_b,y_b)=(0.5,0.3)$ is the center of the bubble, $r=\sqrt{(x-x_b)^2+(y-y_b)^2}$ is the radius, and $\temp_b=0.5$ is the maximum perturbation amplitude. The simulation is performed on a square $\Omega=[0,1]^2$ with $p=4$, $\nu=10^{-6}$, $\alpha=0$ and $\beta=3.411223 \cdot 10^{-3}$. The simulation is performed up to $t=800$ with a time step of $\Delta t=10$. The domain is covered with $\Ni=5616$ triangles. The numerical solution using $R=0,1,2$ is reported in Figure \ref{fig.WV.1}. We can qualitatively observe the same evolution and the same structures as observed in the literature, see \cite{NatConv20} and \cite{Giraldo2008}, with the generation of some small structures when $R=2$, hence showing a less numerical diffusion of the scheme. It is worth to be noted that the first order case $R=0$ is not able to catch the correct evolution of the bubble due to a lack of accuracy in the Lagrangian trajectories. 

\subsection{Density current}
\begin{figure}[!htbp]
	\begin{center}
		\includegraphics[width=0.80\textwidth]{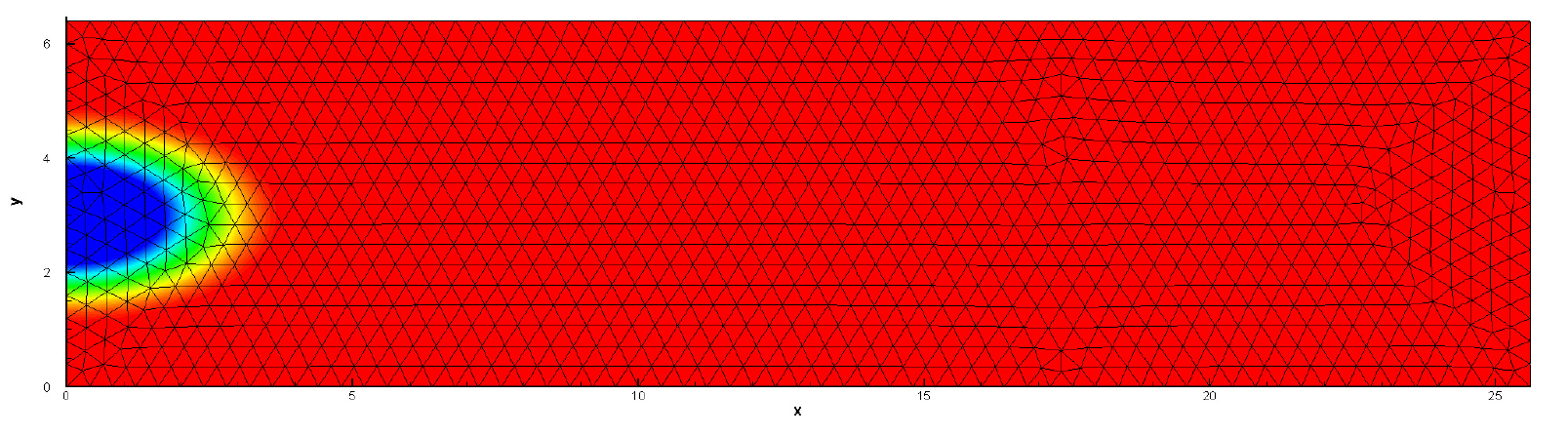}
		\caption{Density current test. Initial condition and computational mesh.}
		\label{fig.DC.1}
	\end{center}
\end{figure}

\begin{figure}[!htbp]
	\begin{center}
		\includegraphics[width=0.32\textwidth]{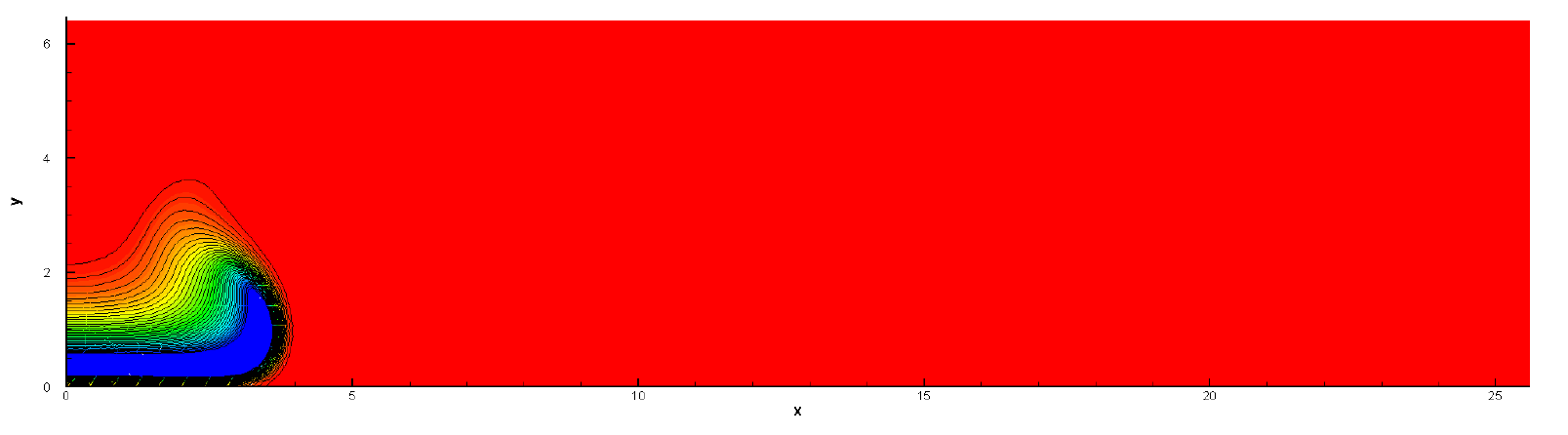}
		\includegraphics[width=0.32\textwidth]{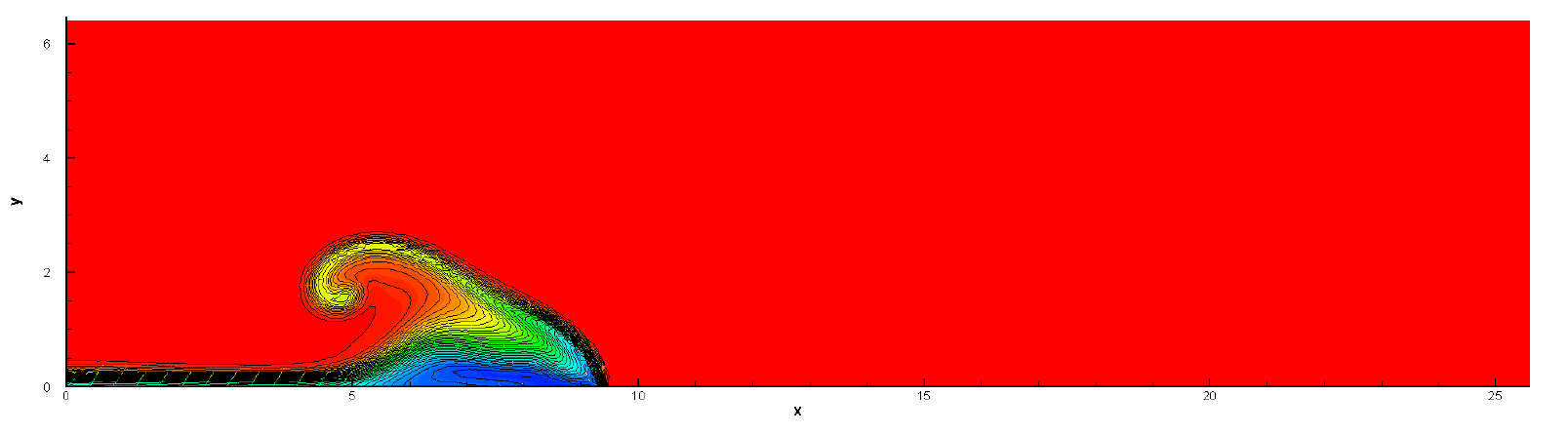}
		\includegraphics[width=0.32\textwidth]{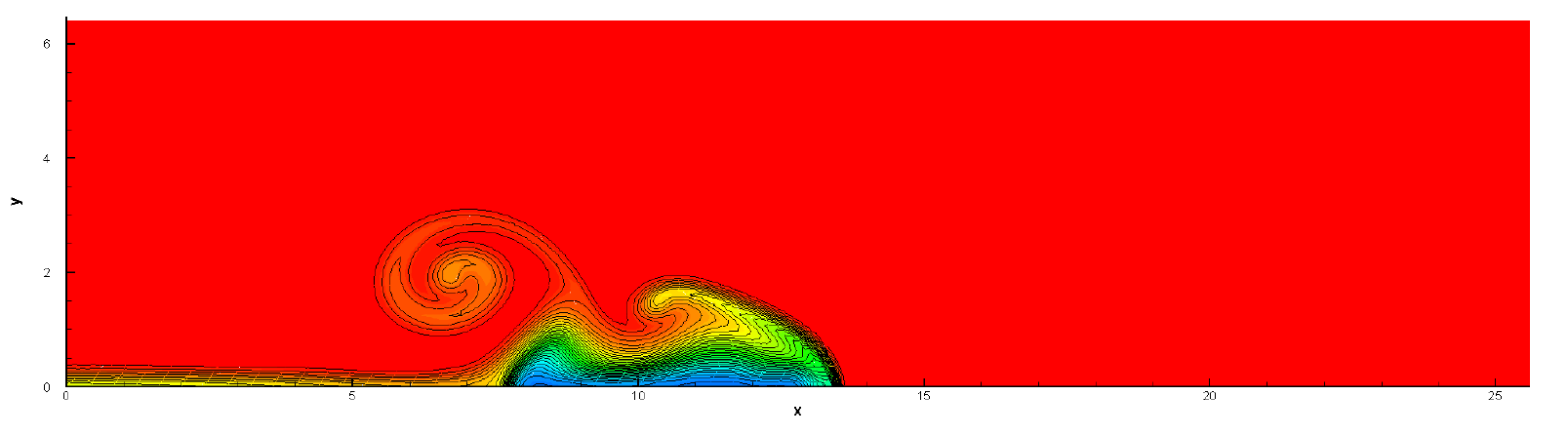} \\
		\includegraphics[width=0.32\textwidth]{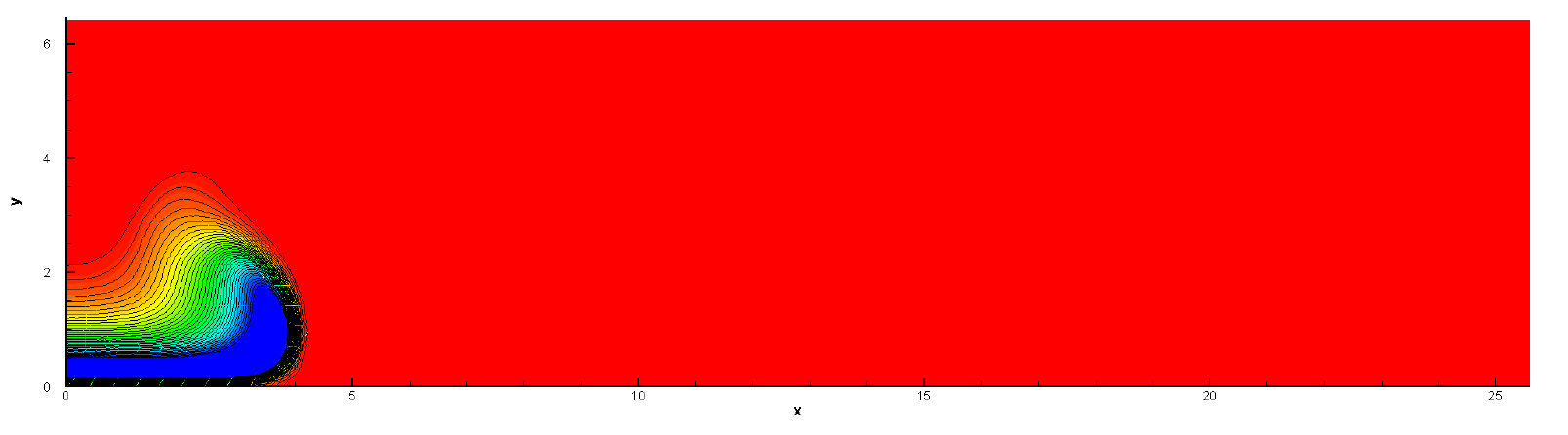}
		\includegraphics[width=0.32\textwidth]{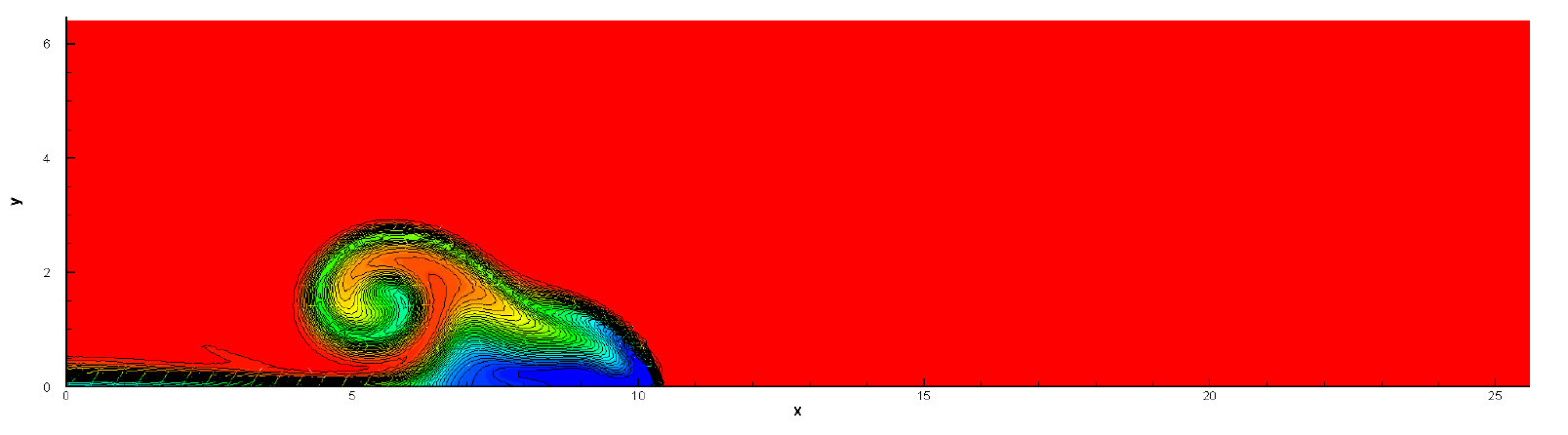}
		\includegraphics[width=0.32\textwidth]{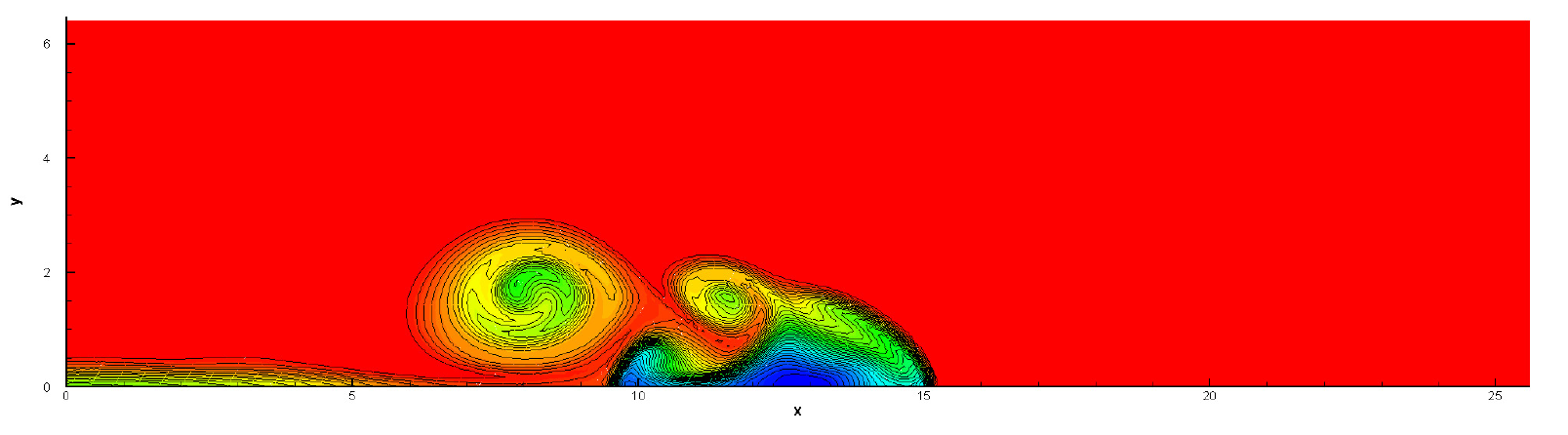} \\
		\includegraphics[width=0.32\textwidth]{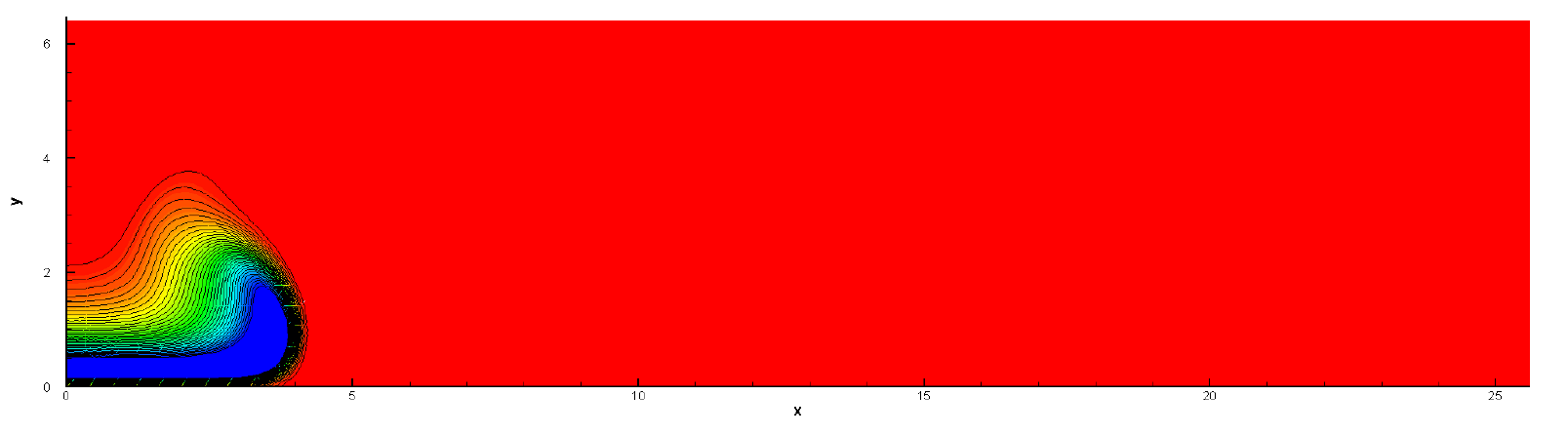}
		\includegraphics[width=0.32\textwidth]{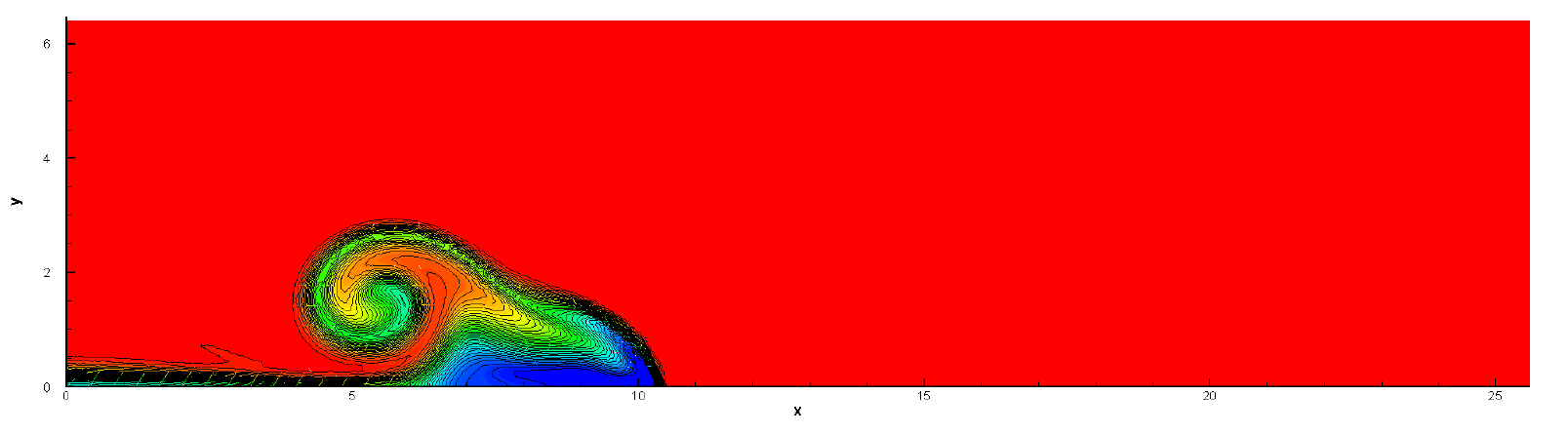}
		\includegraphics[width=0.32\textwidth]{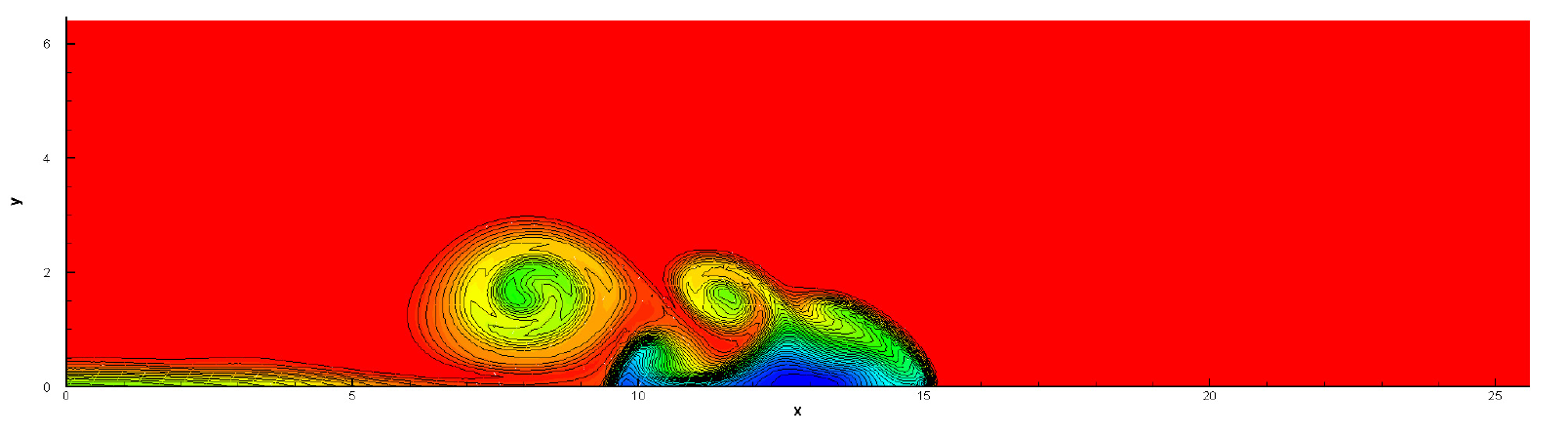} \\
		\caption{Time evolution of the density current problem at times $t=300$ (left column), $t=600$ (central column) and $t=900$ (right column) using $R=0,1,2$ respectively from the first to the third row.}
		\label{fig.DC.2}
	\end{center}
\end{figure}

In order to test the scheme when we have both a density gradient and a viscosity contribution we simulate a density current induced by a cold bubble. This test was originally proposed in \cite{Straka93} and then solved using spectral elements and DG methods in \cite{Giraldo2008}. The temperature gradient is defined as
\begin{eqnarray}
	\temp(x,0)=\left\{
	\begin{array}{lc}
		0 & r>1 \\
		\frac{\temp_b}{2}\left[1 + \cos\left(\pi r\right)\right] & r \leq 1
	\end{array}\right.
\end{eqnarray}
where $\temp_b=-15$ and the radius $r$ is computed as
\begin{eqnarray}
	r=\sqrt{\left(\frac{x-x_c}{x_r}\right)^2+\left(\frac{y-y_c}{y_r}\right)^2},
\end{eqnarray}
with $(x_c,y_c)=(0,3)$ and $(x_r,y_r)=(4,2)$. We finally set $\nu=\alpha=7.5\cdot 10^{-5}$, $\beta=3.41223\cdot 10^{-3}$ and $g=-9.81\cdot 10^{-3}$. The computational domain $\Omega=[0,25.6]\times [0,6.4]$ is covered with $\Ni=2308$ triangles and slip boundary conditions are imposed everywhere. For all the cases we take $p=4$, $t_{end}=900$, $\Delta t =10$ and a different value for $R$, hence analyzing the effect of the time accuracy. The initial condition is plot in Figure \ref{fig.DC.1} while the time evolution for $R=0,1,2$ is depicted in Figure \ref{fig.DC.2} at times $t=300, 600$ and $900$.
The obtained solution is very similar to the one reported in \cite{Giraldo2008} and \cite{Straka93} when we use $R=1$ and $R=2$. In both cases, we can observe the generation of two main circulations past the cold front. In the case $R=0$, only one main circulation is obtained. Furthermore, the front propagates slowly with respect to $R=1$ and $R=2$ and it results to be located far from the value of $15.2$ found in \cite{Straka93}. The main circulation is well resolved for the case $R=1$ and $R=2$ while it results smoothed for the case $R=0$ due to excessive numerical dissipation of first order schemes.

\section{Conclusions}\label{sec:conclusions}
In this paper, we have presented a new class of high order discontinuous-Galerkin (DG) schemes with IMEX time integrators and a semi-Lagrangian (SL) discretization for the advection/transport contribution on two-dimensional unstructured staggered meshes. The DG space discretization is based on a less regular definition of the basis functions on the dual mesh that allows to derive all the differential operators in terms of a limited set of reference operators and simple geometric quantities, with a remarkable saving in terms of memory. The IMEX scheme is developed in the context of DG staggered meshes, which is moreover combined with a semi-Lagrangian technique to account for the convective contribution. According to \cite{TBP22}, the trajectory equation is integrated in the IMEX process to achieve high order of accuracy. We test our SL IMEX-DG schemes for the advection-diffusion, the incompressible Navier-Stokes and the natural convection model. Since the IMEX algorithm belongs to the class of Method Of Lines (MOL) time integrators, the space discretization is decoupled from the temporal one. As such, the resulting implicit systems for the pressure and the viscosity are symmetric and at least semi-positive definite, hence third order in time can be obtained without scarifying those properties, contrarily to the case of fully space-time DG discretizations \cite{TD15}. 

The compactness of the DG operators devised in this work would allow to apply this methodology to moving meshes and $ph$-adaptive grids. Both are interesting applications that will be investigated in future works. Even if other PDE systems, such as the shallow water model, can easily be embedded in this framework, the proposed approach makes use of the non-conservative version of the semi-Lagrangian scheme. In the future, we plan to investigate the use of conservative semi-Lagrangian schemes with staggered DG methods along the lines of \cite{TBP22}, hence permitting to solve also the compressible Navier-Stokes as well as the shallow water equations.

\section*{Acknowledgments}
MT and WB are part of the GNCS INdAM group. MT work was partially funded by the INdAM-GNCS project CUP\_E53C22001930001. WB received financial support by Fondazione Cariplo and Fondazione CDP (Italy) under the grant No. 2022-1895 and by the Italian Ministry of University and Research (MUR) with the PRIN Project 2022 No. 2022N9BM3N.

% % % % % % % % % % % % % % % % % % % % % % % % % % % % % %
% % % % % % % % % % % % % % % % % % % % % % % % % % % % % %
%                  References                             %
\bibliographystyle{elsarticle-num}
\bibliography{./mibiblio}

% % % % % % % % % % % % % % % % % % % % % % % % % % % % % %
% % % % % % % % % % % % % % % % % % % % % % % % % % % % % %
%                   Appendix                              %
% % % % % % % % % % % % % % % % % % % % % % % % % % % % % %
% % % % % % % % % % % % % % % % % % % % % % % % % % % % % %

\appendix

\section{Tableaux of semi-implicit IMEX Runge-Kutta schemes} \label{app.imex}
The IMEX coefficients have been derived in \cite{BosFil2016} and they are reported hereafter for completeness. The left tableau refers to the explicit scheme with $(\tilde{A},\tilde{b},\tilde{c})$, while the right tableau describes the implicit scheme with $(A,b,c)$ in the Butcher table \eqref{eqn.butcher}.

\begin{itemize}
	\item $R=0$, Euler, $\Ns=1$, Formal Order $1$
	
	\begin{equation}
		\begin{array}{c|c}
			0 & 0 \\ \hline & 1
		\end{array} \qquad
		\begin{array}{c|c}
			1 & 1 \\ \hline & 1
		\end{array}
		\label{eqn.IMEX1}
	\end{equation}
	
	\item $R=1$, Diagonally Implicit Runge-Kutta 2, $\Ns=2$, Formal Order $2$ 
	
	\begin{equation}
		\begin{array}{c|cc}
			0 & 0 & 0 \\ \beta & \beta & 0 \\ \hline & 1-\gamma & \gamma
		\end{array} \qquad
		\begin{array}{c|cc}
			\gamma & \gamma & 0 \\ 1 & 1-\gamma & \gamma \\ \hline & 1-\gamma & \gamma
		\end{array}
		\label{eqn.IMEX2}
	\end{equation}
	where $\gamma=1-1/\sqrt{2}$, \hspace{0.05cm} $\beta=1/(2\gamma)$
	
	\item $R=2$, Stiffly accurate Diagonally Implicit Runge-Kutta, $\Ns=4$, Formal Order $3$ 
	\begin{equation}
		\begin{array}{c|cccc}
			0 & 0 & 0 & 0 & 0 \\ \gamma & \gamma & 0 & 0 & 0 \\ 0.717933 & 1.437745 & -0.719812 & 0 & 0 \\ 1 & 0.916993 & 1/2 & -0.416993 & 0 \\ \hline  & 0 & 1.208496 & -0.644363 & \gamma
		\end{array} \qquad
		\begin{array}{c|cccc}
			\gamma & \gamma & 0 & 0 & 0 \\ \gamma & 0 & \gamma & 0 & 0 \\ 0.717933 & 0 & 0.282066 & \gamma & 0  \\ 1 & 0 & 1.208496 & -0.644363 & \gamma \\ \hline  & 0 & 1.208496 & -0.644363 & \gamma
		\end{array}
		\label{eqn.IMEX3}
	\end{equation}
	with $\gamma=0.435866$.

\end{itemize} 

%%%%%%%%%%%%%%%%%%%%%%%%%%%%%%%%%%%%%%%%%%%%%%%%%%%%%%%%
%%%%%%%%%%%%%%%%%%%%%%%%%%%%%%%%%%%%%%%%%%%%%%%%%%%%%%%%
\end{document}